\documentclass[11pt]{amsart}
\usepackage{amsxtra}
\usepackage{amssymb,mathrsfs}
\usepackage{mathtools}
\usepackage{color}
\usepackage{tikz}
\theoremstyle{plain}
\newcommand{\cleqn}{\setcounter{equation}{0}}
\newcommand{\clth}{\setcounter{theorem}{0}}
\newcommand {\sectionnew}[1]{\section{#1}\cleqn\clth}
\newcommand{\nn}{\hfill\nonumber}
\newtheorem{theorem}{Theorem}[section]
\newtheorem{lemma}[theorem]{Lemma}
\newtheorem{definition-theorem}[theorem]{Definition-Theorem}
\newtheorem{proposition}[theorem]{Proposition}
\newtheorem{corollary}[theorem]{Corollary}
\newtheorem{definition}[theorem]{Definition}
\newtheorem{example}[theorem]{Example}
\newtheorem{remark}[theorem]{Remark}
\newtheorem{conjecture}[theorem]{Conjecture}

\newcommand \bth[1] { \begin{theorem}\label{t#1} }
\newcommand \ble[1] { \begin{lemma}\label{l#1} }

\newcommand \bpr[1] { \begin{proposition}\label{p#1} }
\newcommand \bco[1] { \begin{corollary}\label{c#1} }
\newcommand \bde[1] { \begin{definition}\label{d#1}\rm }
\newcommand \bex[1] { \begin{example}\label{e#1}\rm }
\newcommand \bre[1] { \begin{remark}\label{r#1}\rm }
\newcommand \bcj[1] { \begin{conjecture}\label{j#1}\rm }

\renewcommand {\eth} { \end{theorem} }
\newcommand {\ele} { \end{lemma} }

\newcommand {\epr} { \end{proposition} }
\newcommand {\eco} { \end{corollary} }
\newcommand {\ede} { \end{definition} }
\newcommand {\eex} { \end{example} }
\newcommand {\ere} { \end{remark} }
\newcommand {\ecj} { \end{conjecture} }

\newcommand {\enota} { \end{notation} }
\newcommand \thref[1]{Theorem \ref{t#1}}
\newcommand \leref[1]{Lemma \ref{l#1}}
\newcommand \prref[1]{Proposition \ref{p#1}}
\newcommand \coref[1]{Corollary \ref{c#1}}

\newcommand \deref[1]{Definition \ref{d#1}}

\newcommand \lb[1]{\label{#1}}
\def \Cset {{\mathbb C}}

\def \Zset {{\mathbb Z}}

\def \Abb {{\mathcal A}_q^{1/2}}

\def \A {{\boldsymbol{\mathsf A}}} 
\def \C {{\boldsymbol{\mathsf C}}} 
\def \U {{\boldsymbol{\mathsf U}}} 

\def \B  {{\widetilde{B}}}
\newcommand \ex {{\bf{ex}}}
\newcommand \inv {{\bf{inv}}}
\newcommand \ninv {{\bf{ninv}}}
\newcommand \nc {{\bf{nc}}}

\def \x {{\wt{\mathbf{x}}}} %seed notation

\def \AA  {{\mathcal{A}}}           %mathcal

\def\DD{{\mathcal{D}}}

\def \VV {{\mathcal{V}}}

\def \SS {{\mathcal{S}}}
\def \LL {{\mathcal{L}}}
\def \Lsc {{\mathscr{L}}}
\def \TT {{\mathcal{T}}} 
\def \ZZ {{\mathcal{Z}}}

\def \Osc {{\mathscr{O}}}

\def \FF {{\mathcal{F}}}

\def \ZZ {{\mathcal{Z}}}

\def \vpi {\varpi}
\def \la {\lambda}
\def \La {\Lambda}

\def \Om {\Omega}

\def \Ga {\Gamma}

\def \ep {\varepsilon}
\def \mt  {\mapsto}
\def \hra {\hookrightarrow}

\def \rcor {\rangle}
\def \lcor {\langle}
\def \ol {\overline}
\def \wt {\widetilde}

\def \id { {\mathrm{id}} }

\def \reg   { {\mathrm{reg}} }
\def \sing { {\mathrm{sing}}}
\def \n  {\mathfrak{n}}
\def \mm  {\mathfrak{m}}

\DeclareMathOperator \Cl  {{\mathrm{Cl} }}
\DeclareMathOperator \Span { {\mathrm{Span}} }

\DeclareMathOperator \diag { {\mathrm{diag}} }
\DeclareMathOperator \Ker { {\mathrm{Ker}} }

\DeclareMathOperator \rk { {\mathrm{rk}} }

\DeclareMathOperator \MaxSpec { {\mathrm{MaxSpec}} }

\DeclareMathOperator \Supp { {\mathrm{Supp}} }
\DeclareMathOperator \Der { {\mathrm{Der}}}
\renewcommand \Im { {\mathrm{Im}} }

\renewcommand \max { {\mathrm{max}} }

\newcommand \Fract {{\mathrm{Fract}}}
\begin{document}
\title
[Poisson geometry and Azumaya loci of cluster algebras]
{Poisson geometry and Azumaya loci \\ of cluster algebras}
\author[G. Muller]{Greg Muller}
\thanks{B.N. has been supported by an AMS-Simons travel grant and M.Y.  by NSF grants DMS-2131243 and DMS--2200762.}
\address{Department of Mathematics, University of Oklahoma, Norman, OK 73019}
\email{gmuller@ou.edu}
\author[B. Nguyen]{Bach Nguyen}
\address{
Department of Mathematics, Xavier University of Louisiana, New Orleans, LA 70125}
\email{bnguye22@xula.edu}
\author[K. Trampel]{Kurt Trampel}
\address{
Department of Mathematics, University of Notre Dame, Notre Dame, IN, 46556}
\email{ktrampel@nd.edu}
\author[M. Yakimov]{Milen Yakimov}
\address{
Department of Mathematics, Northeastern University, Boston, MA 02115}
\email{m.yakimov@northeastern.edu}
\begin{abstract} There are two main types of objects in the theory of cluster algebras: the upper cluster algebras $\U$
with their Gekhtman--Shapiro--Vainshtein Poisson brackets and their root of unity quantizations $\U_\ep$. 
On the Poisson side, we prove that (without any assumptions) the spectrum of every finitely generated upper cluster algebra $\U$ with its GSV Poisson structure 
always has a Zariski open orbit of symplectic leaves and give an explicit description of it. 
On the quantum side, we describe the fully Azumaya loci of the quantizations $\U_\ep$ under the assumption that 
$\A_\ep = \U_\ep$ and $\U_\ep$ is a finitely generated algebra. All results allow frozen variables to be either inverted or not.
\end{abstract}
\maketitle
%%%%%%%%%%%%%%%%%%%%%%%%%%%%%%%%%%%%%%%%%%%%%%%%%%%%%%%%%%%%%%%%%%%%%%%%%%%
\sectionnew{Introduction}
\lb{Intro}
\subsection{Setting}
\label{1.1}
Cluster Algebras were defined by Fomin and Zelevinsky \cite{FZ1} in 2001, and since then have played a prominent role in 
many areas of mathematics and mathematical physics. There are two main classes of algebras in this theory. The algebras in these classes are defined 
from an integer matrix $\B$ (called {\em{exchange matrix}}) of size $N \times \ex$ where $N$ is a positive integer and $\ex \subset [1,N]$ 
is a set of mutable (exchangeable) variables. Furthermore, one partitions
\[
[1,N] \backslash \ex = \inv \sqcup \ninv.
\]
The two subsets $\inv$ and $\ninv$ will index the inverted and non-inverted frozen variables. It is important to allow this degree of flexibility because 
many key examples in Lie theory require that not all frozen variables are inverted.  
\begin{enumerate}
\item The {\em{upper cluster algebra}} $U(\B,\inv)$ is defined as the intersection
\[
\U(\B,\inv)= \bigcap_{(\x', \B') \sim (\x, \B)} \Cset[(x'_k)^{\pm 1}, x'_i ; k \in \ex \sqcup \inv, i \in \ninv], 
\]
where the intersection ranges over all seeds $(\x':= (x'_1, \ldots, x'_N), \B)$ in the mutation class of the initial seed
$(\x=(x_1, \ldots, x_N), \B)$. Gekhtman, Shapiro and Vainshtein proved in \cite{GSV} that, 
if the exchange matrix $\B$ is compatible with a skew-symmetric integer matrix $\La$ of size $N \times N$ (cf. \eqref{compatibility}),
then $\U(\B,\inv)$ admits a canonical {\em{Poisson structure}}. The cluster 
variables in each seed are in log-canonical form, meaning that
\[
\{ x'_k, x'_i \} = \la'_{ki} x'_k x'_i, \quad \forall k,i \in [1,N]
\]
for some skew-symmetric integer matrix $(\la'_{ki})$ depending on the seed.
\item Let $\ep^{1/2} \in \Cset$ be a primitive $\ell$-th root of unity for a positive integer $\ell$. 
The root of unity upper quantum cluster algebra $\U_\ep(M_\ep,\B,\inv)$ is a non-commutative algebra, 
defined in a similar way to $\U(\B,\inv)$, by intersecting 
mixed quantum tori/quantum affine spaces 
\[
\U_\ep(M_\ep, \B, \inv):= \bigcap_{(M'_\ep,\B')\sim (M_\ep,\B)} \TT_\ep(M'_\ep)_\geq,
\]
where $M'_\ep$ are root of unity toric frames, see Sect. \ref{3.1}. Each such algebra has a 
canonical central subalgebra $\C\U_\ep(M_\ep, \B, \inv)$ obtained by intersecting the mixed Laurent polynomial 
rings in the $\ell$-th powers of the generators of  $M'_\ep$. Under mild assumptions on $\ell$,  
\begin{equation}
\label{u-CUep-isom}
\C\U_\ep(M_\ep, \B, \inv) \cong \U(\B, \inv),
\end{equation} 
see Sect. \ref{3.1}. 
\end{enumerate}
There is a canonical action of the torus
\[
T(\B) := ( \Cset^\times)^{\dim \Ker \B}
\] 
on both algebras $\U(\B, \inv)$ and $\U_\ep(M_\ep, \B, \inv)$ by algebra automorphisms. In the first case, $T(\B)$ acts by 
Poisson automorphisms. In the second case, it preserves the central 
subalgebra $\C\U_\ep(M_\ep, \B, \inv)$ and the isomorphism \eqref{u-CUep-isom} intertwines the two actions. 
%%%%%%%%%%%
\subsection{Results on the Poisson geometric side}
\label{1.2}
The affine Poisson variety corresponding to the GSV Poisson structure
\[
Y(\B):= \MaxSpec \U(\B,\inv)
\]
is of much interest in Lie theory and integrable systems (we suppress the dependance of $Y(\B)$ on $\inv$ for brevity).
However, little is known about its global geometry. One only knows that the cluster tori inside $Y(\B)$ are regular Poisson \cite{GSV}. 
This is a local result because the symplectic leaves inside each cluster torus never entirely belong to the cluster torus. So, 
general Hamiltonian flows on $Y(\B)$ leave each cluster torus after some time. 

In Lie theory and combinatorics, one knows \cite{GY0} that each Schubert cell $X_w^\circ$ in the full flag variety $G/B_+$ 
of a complex simple Lie group $G$, equipped with the standard Poisson structure, has a dense torus orbit of symplectic leaves which equals the complement of the 
Richardson divisor \cite{KLS} of $X_w^\circ$. The latter equals the union of the closures in $X_w^\circ$ of the 
open Richardson varieties $R_{w,s_i}$ for $i$ ranging over the support of the Weyl group element, see Sect. \ref{7.1}.

Remarkably, such a fact holds for the spectrum of every upper cluster algebra $\U(\B, \inv)$ with the GSV Poisson structure,
without any assumptions on $\U(\B, \inv)$ except for finite generation which is needed to be even able to talk about symplectic leaves.
This is the topic of our first main result on the the global description of the $T(\B)$-orbit of the symplectic leaves $\SS$ of $(Y(\B), \pi)$ of maximal dimension,
\[
T(\B) \cdot \SS.
\] 

\noindent
{\bf{Theorem A.}}
{\em{Assume that $\U(\B, \inv)$ is a finitely generated upper cluster algebra for which the exchange matrix $\B$ admits a compatible skew-symmetric integer matrix. 
Then the affine Poisson variety $(Y(\B), \pi)$ has a Zariski open $T(\B)$-orbit of symplectic leaves, which equals
\[
Y(\B)^\reg \backslash \bigcup_{i \in \ninv} \VV (\x_i),
\]
where $Y(\B)^\reg$ is the nonsingular part of $Y(\B)$. 
}}
\medskip

%Here the symplectic leaves of a singular Poisson variety $Y$ are defined using Polishchuk's method \cite{P} that the singular part of $Y^\sing$ 
%is a Poisson subvariety of $Y$, and thus the leaves are constructed as the leaves of $Y^\reg, (Y^\sing)^\reg, \ldots$

The theorem is proved by using sections of the anticanonical bundle of $Y(\B)^\reg$ coming from the GSV Poisson structure $\pi$ and the action of 
$T(\B)$, coupled with normality of $Y(\B)$. 

At this point it might be tempting to conjecture that a much stronger result than Theorem A holds, namely that there are 
finitely many $T(\B)$-orbits of symplectic leaves of $Y(\B)$ (or $Y(\B)^\reg$). However, that is not correct. 
For large classes of Belavin--Drinfeld Poisson structures \cite[Sect. 3.2]{CP} on ${\mathrm{GL}}_n(\Cset)$, 
Gekhtman, Shapiro and Vainshtein proved \cite{GSV1,GSV2} that the coordinate ring $\Cset[{\mathrm{GL}}_n]$ admits 
upper cluster algebra structures with compatible Poisson structures given by the ones in the list. 
Their symplectic leaves were classified in \cite{Y}, where it was proved that they are classified 
by Weyl group datum and twisted conjugacy classes of reductive groups. Because of the last bit of data, in general,
those Poisson structures have infinitely many torus orbits of symplectic leaves.
%%%%%%%%%%%
\subsection{Results on the quantum side}
\label{1.3}
It was proved in \cite[Theorem B]{HLY} that, if the $\Cset$-algebra $\U_\ep(M_\ep, \B, \inv)$ is finitely generated, 
then it is a finitely generated module over $\C\U_\ep(M_\ep, \B, \inv)$ and $\C\U_\ep(M_\ep, \B, \inv) \cong \U(\B, \inv)$ is a finitely 
generated commutative $\Cset$-algebra. The root of unity upper quantum cluster algebras $\U_\ep(M_\ep, \B, \inv)$ form 
a vast family of algebras that includes as special cases many important classes of quantum algebras at roots of unity 
arising in Lie theory and topology. A fundamental open problem for them is to classify their irreducible representations. 
By \cite[Theorem III.1.6]{BrG} all of their representations are finite dimensional of dimension less than or equal to the PI degree 
of $\U_\ep(M_\ep, \B, \inv)$. In this generality, the problem of classifying the irreducible representations of $\U_\ep(M_\ep, \B, \inv)$ is extremely difficult. 
The first step towards its resolution 
is to classify the representations of maximal dimension. For this, Brown and Gordon defined \cite[Sect. 1.3]{BG1} the  
{\em{fully Azumaya locus}} of a finitely generated prime algebra $R$ with respect to a central subalgebra $Z$ such that 
$R$ is a finitely generated $Z$-module. This locus is a Zariski open subset of $\MaxSpec Z$ consisting of those maximal ideals $\mm$ of $Z$ 
such that all irreducible representations of $R$ annihilated by $\mm$ have maximal dimension (equal to the PI degree of $R$), 
see Sect. \ref{2.1} for details. 

Our second main result gives an explicit description of the fully Azumaya locus of the root of unity upper 
quantum cluster algebras $\U_\ep(M_\ep, \B, \inv)$ with respect to their central subalgebras $\C\U_\ep(M_\ep, \B, \inv) \cong \U(\B, \inv)$:
\medskip

\noindent
{\bf{Theorem B.}} {\em{Assume that $\U_\ep(M_\ep, \B, \inv)$ is a finitely generated strict root of unity upper quantum cluster algebra such that 
the order of $\ep^{1/2}$ is odd and coprime to the diagonal entries of the skew-symmetrizing matrix $D$ 
for the principal part of the exchange matrix $\B$ and 
\begin{equation}
\label{A=U}
\A_\ep(M_\ep, \B,\inv) = \U_\ep(M_\ep, \B,\inv).
\end{equation}
Then the fully Azumaya locus $\AA$ of 
$\U_\ep(M_\ep, \B, \inv)$ with respect to the central subalgebra
\[
\C\U_\ep(M_\ep, \B, \inv) \cong \U(\B, \inv) 
\]
satisfies
\[
Y(\B)^\reg \backslash \bigcup_{i \in \ninv} \VV (\x_i) \subseteq \AA \subseteq 
Y(\B) \backslash \bigcup_{i \in \nc} \VV (\x_i),
\]
where $\nc$ denotes the set of those non-inverted frozen variables $M(e_i)$ that are not in the center on $\U_\ep(M_\ep, \B, \inv)$.
}}
\medskip

The assumption that $\U_\ep(M_\ep, \B, \inv)$ is a finitely $\Cset$-algebra is needed to be even able to define the fully Azumaya locus.
The strictness assumption on $\U_\ep(M_\ep, \B,\inv)$ means that $\B$ admits a skew-symmetric integer matrix which is compatible 
with $\B$  over $\Zset$ and not just over $\Zset/\ell$, which is what is needed to 
define root of unity quantum cluster algebras in general. This assumption ensures that there is an associated upper quantum cluster algebra 
$\U_q(M_q, \B,\inv)$ in the sense of Berenstein and Zelevinsky \cite{BerZe}. The assumption \eqref{A=U} is needed to ensure that 
$\U_\ep(M_\ep, \B,\inv)$ is a specialization of $\U_q(M_q, \B,\inv)$. This in turn is used to 
construct a Poisson order structure on the pair $(\U_\ep(M_\ep, \B,\inv), \C\U_\ep(M_\ep, \B,\inv))$ in the sense of 
Brown and Gordon \cite{BG1} to be able to link Poisson geometry to representation theory. A property of the form $\A=\U$ 
was established in numerous situations on the classical and quantum levels \cite{CGGLSS,GY1,GY1,GY2,GY3,M1,M2,SW}. 
Undoubtably, these methods will be extended in the future to show that the assumption \eqref{A=U} is satisfied in broad generality.

In \cite{BG1a} the fully Azumaya loci of all quantum function algebras of complex simple Lie groups $G$ at roots of unity were determined. 
In \cite{GJS} a second proof of the Bonahon--Wang unicity conjecture \cite{BW} was given, stating that the Azumaya loci of the skein algebras of oriented surfaces 
at roots of unity contain the smooth parts of the spectra of their centers (this was first proved in \cite{FKBL}). 
Both \cite{BG1a,GJS} relied on Poisson orders. 
The result of \cite{BG1a} is that the fully Azumaya locus inside $G$ is the open double Bruhat cell (a torus orbit of symplectic leaves), while the proofs of \cite{GJS} relied on  
an argument that in that situation, the smooth part in question is a single symplectic leaf. Our Theorems A and B prove that such phenomena hold in much wider 
generality in the world of cluster algebras, rather than quantum groups or skein algebras.
The zero loci of frozen variables and singular part of $Y(\B)$ were exactly the points that were thrown out 
in \cite{BG1a} and \cite{GJS}.
\medskip

\noindent 
{\bf{1.4 Organization of the paper and notation.}} The paper is organized as follows. Section 2 provides background on Azumya loci, 
Poisson orders, cluster algebras, compatible Poisson structures and quantum cluster algebras. Section 3 contains background on 
root of unity quantum cluster algebras and results on torus actions on them and the related Poisson cluster algebras.
Section 4 proves Theorem A. Section 5 construct Poisson orders on root of unity quantum cluster algebras. 
Section 6 proves Theorem B. Section 7 discusses the motivation for the main theorem from the stand point of 
Richardson varieties, the special case of the theorems on acyclic cluster algebras, and a Kronecker type 
example.

The following notation will be used throughout the paper. The standard basis of $\Zset^N$ will be denoted by 
\[
e_1, \ldots, e_N.
\]
The dot product on $\Zset^N$ will be denoted by $\mu \cdot \nu$. The transpose of a matrix $B$ will be denoted by $B^\top$.

For a complex affine algebraic variety $Y$, its singular part will be denoted by $Y^\sing$ and its smooth part by $Y^\reg$. Given a regular function 
$f \in \Cset[Y]$, we will denote by $\VV(f)$ the zero locus of $f$. 
%%%%%%%%%%%
\sectionnew{Preliminaries on Poisson orders and cluster algebras}
\lb{prelim}
In this section we gather background material on Poisson orders and cluster algebras that 
will be used in the paper.
\subsection{Poisson orders}
\label{2.1}
We follow Brown and Gordon \cite{BG2}.
Let $R$ be a $\Cset$-algebra and $Z$ be a central subalgebra of $R$.
\bde{poi order} \cite{BG2}
The pair $(R,Z)$ is called a {\em{Poisson order}} if $R$ is a finitely generated $Z$-module satisfying the following conditions:
\begin{enumerate}
\item[(i)] $Z$ is has the structure of a Poisson algebra with bracket $\{\cdot,\cdot\}$;
\item[(ii)] There exists a $\Cset$-linear map $\partial: Z\longrightarrow \Der_\Cset(R)$ such that $\partial_z|_Z=\{z,-\}$ for all $z\in Z$.
\end{enumerate}
\ede
By assumption (ii), the Poisson structure is uniquely determined from the linear map $\partial$. 
Because of this we will denote Poisson orders as triples $(R,Z, \partial)$. 

{\em{Restriction of Poisson orders:}} For a Poisson order $(R, Z, \partial)$, if $C$ is a Poisson subalgebra of $Z$ with respect to 
the underlying Poisson structure, then $(R,C, \partial|_C)$ is also a Poisson order. 

{\em{Poisson orders from specialization:}} The following is a well known fact for obtaining Poisson order structures from specialization,
see e.g. \cite[Sect. 2.2]{BG2}:
\ble{Pord-spec} 
Assume that $R$ and $S$ are $\Cset$-algebras and $\eta : S \to R$ is a surjective $\Cset$-algebra homomorphism with kernel $(h)=hS$ 
for a regular central element $h \in S$. Choose a $\Cset$-linear map
\[
\iota : \ZZ(R) \to S
\]
such that $\eta \circ \iota = \id_{\ZZ(R)}$. If $R$ is a finitely generated $\ZZ(R)$-module, then 
the pair $(R, \ZZ(R))$ admits a Poisson order structure with $\partial_z : \ZZ(R) \to \Der_\Cset(R)$ 
given by
\begin{equation}
\label{partial-spec}
\partial_z (r) := \eta \Big( \frac{\iota(z) \wt{r} - \wt{r} \iota(z)}{h} \Big), \quad \forall z \in \ZZ(R), r \in R,
\end{equation}
where $\wt{r}$ is any preimage of $r$ inder $\eta$. Its underlying Poisson structure is given by 
\[ 
\{z_1, z_2 \} := \eta \Big( \frac{\iota(z_1) \iota(z_2) - \iota(z_2) \iota(z_1)}{h} \Big), \quad \forall z_1, z_2 \in \ZZ(R).
\]
This Poisson structure is independent on the choice of $\Cset$-linear section $\iota : \ZZ(R) \to S$.
\ele
In \eqref{partial-spec}, $\iota(z) \wt{r} - \wt{r} \iota(z) \in (h)$ because $z \in \ZZ(R)$ and $\eta$ is an 
algebra homomorphism. The right hand side of  \eqref{partial-spec} is independent on the choice of preimage $\wt{r}$
by a similar argument.

Recall that an affine Poisson variety is an affine variety $X$ whose coordinate ring is equipped with a Poisson algebra structure. 
Its singular part $X^\sing$ is automatically an affine Poisson variety as well \cite[Corollary 2.4]{P}. The symplectic leaves of 
a complex affine Poisson variety $X$ are defined recursively as the symplectic leaves of the smooth complex manifold $X^\reg$ 
together with the symplectic leaves of the lower dimensional Poisson variety $X^\sing$.

\bth{P-ord-symp-leaves} {\em{(Brown--Gordon)}} \cite{BG2} Assume that $(R,Z)$ is a complex Poisson order and $\mm, \mm' \in \MaxSpec Z$ lie in the 
same symplectic leaf. Then we have the isomorphism of finite dimensional complex algebras 
\[
R/\mm R \cong R/ \mm' R.
\] 
\eth

The theorem proved in \cite{BG2}  has a stronger conclusion using Poisson cores, but we will not need that fact in this paper. 

For the rest of this subsection we assume that 

(*) {\em{$R$ is a finitely generated prime $\Cset$-algebra, which is a finitely generated $Z$-module over a central subalgebra $Z$.}} 

By the Artin--Tate Lemma (see e.g. \cite[Sect. I.13.4]{BrG}), $Z$ is a finitely generated $\Cset$-algebra.  
The primeness assumption on $R$ implies that $\ZZ(R)$, and thus also $Z$, are integral domains. Hence, $\MaxSpec \ZZ(R)$ and $\MaxSpec Z$ 
are irreducible affine varieties. 
 
Recall that $\mm \in \MaxSpec \ZZ(R)$ is in the {\em{Azumaya locus of $R$}} if $R_\mm$ is an Azumaya algebra over 
$Z_\mm$. This is equivalent to saying that $R$ has an irreducible module, annihilated by $\mm$, of maximal dimension 
among the irreducible $R$-modules (which equals the PI degree of $R$); 
such a representation is automatically unique (see \cite[Theorem III.1.6]{BrG}). 

We have the canonical map
\begin{equation}
\label{ZRtoZ}
\MaxSpec \ZZ(R) \to \MaxSpec Z,
\end{equation}
induced by the inclusion $Z \subseteq \ZZ(R)$. 
\bde{fullyAzumaya}
A point $\mm \in \MaxSpec Z$ is said to be in the {\em{fully Azumaya locus of $R$ with respect to $Z$}} if all of its preimages are 
in the Azumaya locus of $R$, see \cite[Sect. 1.3]{BG1}. In other words, one requires that all irreducible modules of $R/ \mm R$ have 
dimensions equal to the PI degree of $R$.
\ede
The map \eqref{ZRtoZ} is closed by \cite[Lemma III.1.5]{BrG}. 
This and the fact that the Azumaya locus of $R$ is open, and hence dense in $\MaxSpec \ZZ(R)$, imply 

\ble{FAzumaya-dense} In the above setting, the fully Azumya locus of $R$ with respect to $Z$ is an open and hence dense subset of $\MaxSpec Z$.
\ele
This fact is stated in \cite[Proposition III.4.10]{BrG} under a Hopf algebra assumption on $R$ and $Z$, but this assumption is not used in its proof. 

Finally, we have the following corollary of \thref{P-ord-symp-leaves}:
\bco{P-ort-FAzumaya} 
Assume that $R$ is a finitely generated $\Cset$-algebra and
$(R,Z)$ is a complex Poisson order. Then the fully Azumaya locus of $R$ with respect to $Z$ is a union of symplectic leaves of $\MaxSpec Z$.
\eco
\begin{proof} If $\mm, \mm' \in \MaxSpec Z$ lie in the same symplectic leaf of $\MaxSpec Z$, 
then $R/\mm R \cong R/ \mm' R$ by \thref{P-ord-symp-leaves}. In particular, if all irreducible representations of the 
algebra $R/ \mm R$ have dimensions equal to the PI degree of $R$, then the same holds for the algebra 
$R/ \mm' R$.
\end{proof}
\subsection{Cluster algebras of geometric type}
\label{2.2}
We follow Berenstein, Fomin and Zelevinsky \cite{FZ1,BFZ05} with the exception that all algebras are defined over $\Cset$ instead of $\Zset$. 
Let $N$ be a positive integer, $\ex\subseteq [1,N]$ and $\FF$ be a purely transcendental field extension of $\Cset$ 
of transcendence degree $N$. We say that a pair $(\x, \B)$ is a {\em{seed}} if 
\begin{enumerate}
\item[(i)] $\x = \{x_1, \dots, x_N\}$ is a transcendence basis of $\FF$ over $\Cset$;
\item[(ii)] $\B\in M_{N\times \ex}(\Zset)$ and its $\ex\times \ex$ submatrix $B$, called the principal part of $\B$, 
is skew-symmetrizable for some matrix $D=\diag(d_j,~j\in\ex)$ with diagonal entries $d_j \in \Zset_+$. %, which will be assumed to be relatively prime. 
\end{enumerate}
The elements $x_i \in \FF$ are called {\em{cluster variables}}. A matrix $\B$ satisfying the condition in (ii) is called an {\em{exchange matrix}}.
For each $k\in \ex$, the mutation of $\B$ in the direction of $k$ is the matrix $\mu_k(\B)$, where

\[\mu_k(\B) = (b_{ij}'):=
\begin{cases}
-b_{ij} & \text{if } i = k \text{ or } j=k \\
b_{ij} + \frac{|b_{ik}|b_{kj} + b_{ik}|b_{kj}|}{2} & \text{otherwise.}
\end{cases}
 \]
For a choice of sign, $s = \pm$, define the matrices $E_s \in M_N(\Zset)$ and $F_s \in M_\ex(\Zset)$ to be
\[
E_s := (e_{ij}) =
\begin{cases}
\delta_{ij} &  \text{if } j \neq k \\
-1 & \text{if } i = j = k \\
\max(0,-sb_{ik}) & \text{if } i \neq j = k,
\end{cases} \; \; 
F_s := (f_{ij})=
\begin{cases}
\delta_{ij} &  \text{ if } i \neq k \\
-1 & \text{ if } i = j = k \\
\max(0,sb_{kj}) & \text{ if } j \neq i = k.
\end{cases}
\]
Then we also have $\mu_k(\B)= E_s \B F_s$ for both $s = \pm$.

The mutation of a seed $(\x, \B)$ in the direction of $k \in \ex$ 
is defined to be $\mu_k (\x, \B) := (\x', \mu_k(\B))$ where 
\begin{equation}
\label{classic-mut}
\x' := \{x_k'\} \cup \x \backslash \{x_k\} 
\quad
\mbox{and}
\quad
x_k x_k' := \prod_{b_{ik}>0}x_i^{b_{ik}} + \prod_{b_{ik}<0}x_i^{-b_{ik}}.
\end{equation}
The pair $\mu_k (\x, \B)$ is also a seed, the principal part of $\mu_k(\B)$ equals $\mu_k(B)$, and $\mu_k(B)$ is skew-symmetrizable 
with respect to the same matrix $D$ that skew-symmetrizes $B$.

Mutation is an involution. Two seeds are {\em{mutation equivalent}}, $(\x', \B')\sim (\x'', \B'')$, if one can be obtained from the other by a finite sequence of mutations.
Any seed which is mutation-equivalent to $(\x, \B)$ contains $x_i$ for $i \in [1,N]\backslash \ex $ and we call these, {\em{frozen variables}}. 

Fix a decomposition of the set of frozen variables into a disjoint union of two sets: 
\[
[1,N]\backslash \ex= \inv \sqcup \ninv.
\]
The set $\inv$ will index the set of those frozen variables that will be {\em{inverted}}. The set $\ninv$ will index the {\em{non-inverted frozen variables}}.

The {\em{cluster algebra}} $\A(\B, \inv)$ is defined to be the $\Cset$-subalgebra of $\FF$ generated by 
the cluster variables in all the seeds $(\x',\B')\sim( \x, \B)$ together with $\{ x_i^{-1} \mid i \in \inv\}$.
For a seed $(\x',\B')\sim (\x,\B)$ denote the Laurent polynomial ring
\[
\LL(\x') = \Cset[(x'_k)^{\pm 1}; 1 \leq k \leq N]
\]
and its mixed polynomial/Laurent polynomial subring
\begin{equation}
\label{Lx'}
\LL(\x')_\geq = \Cset[(x'_k)^{\pm 1}, x'_i; k \in \ex \sqcup \inv, i \in \ninv].
\end{equation}
The {\em{upper cluster algebra}} $\U(\B, \inv)$ is the intersection 
\begin{equation}
\label{upper-c-a}
\U(\B,\inv)= \bigcap_{(\x', \B') \sim (\x, \B)} \LL(\x')_\geq.
\end{equation}
By the {\em{Laurent phenomenon}}  \cite{FZlp}, we have $\A(\B, \inv) \subseteq ~\U(\B, \inv)$. 

We will need the algebra \eqref{upper-c-a} in the special case when all frozen variables are inverted:
\[
\U(\B):= \U(\B,  [1,N]\backslash \ex) = \bigcap_{(\x', \B') \sim (\x, \B)} \LL(\x').
\]
It is easy to verify that the latter is obtained as a localization:
\[
\U(\B) = \U(\B, \inv)[x_i^{-1}; i \in \ninv].
\]
\subsection{Quantum cluster algebras}
\label{2.3}
We follow Berenstein and Zelevinsky \cite{BerZe} with the exception that we allow for an arbitrary subset of frozen quantum 
cluster variables not to be inverted and we work over the algebra
\[
\Abb:= \Cset[q^{\pm 1/2}]
\]
instead of $\Zset[q^{\pm 1/2}]$.

By abuse of notation, we will identify a skew-symmetric bilinear form $\La: \Zset^N \times \Zset^N \to \Zset$ 
with the skew-symmetric integer matrix with entries $\La(e_i, e_j)$. Recall that $e_1, \ldots, e_N$ denotes the standard basis of $\Zset^N$.

The \emph{based quantum torus} $\TT_{q}(\La)$ associated with $\La$ is the $\Abb$-algebra with the $\Abb$-basis 
$\{ \hspace{1pt} x^f \hspace{1pt} | \hspace{1pt} f \in \Zset^N \hspace{1pt} \}$ and multiplication given by $x^f x^g = q^{\La(f,g)/2} x^{f+g},$ 
where $f,g \in \Zset^N$. 

A \emph{toric frame} $M_q$ for a division $\Cset(q^{1/2})$-algebra $\FF_q$ 
is a map $M_q : \Zset^N \to \FF_q$ for which there exists an $\Abb$-algebra embedding $\phi: \TT_q(\La) \hra \FF_q$
for some skew-symmetric matrix $\La \in M_N(\Zset)$, such that $\phi(x^f)=M_q(f)$ for all $f\in \Zset^N$ and $\FF_q \simeq \Fract\left( \phi( \TT_q(\La) ) \right)$.
Denote by $\La_{M_q}$ the skew-symmetric matrix (bilinear form) of a toric frame $M_q$, defined by
\[
M_q(f) M_q(g) = q^{\La(f,g)/2} M_q(f+g), \quad \forall f, g \in \Zset^N.
\]
For a torus frame $M_q$, the image of $\phi$ in $\FF_q$ will be denoted by $\TT_q(M_q)$; its basis is  $\{ \hspace{1pt} M_q(f) \mid f \in \Zset^N \hspace{1pt} \}$. 
We have the isomorphism of quantum tori $\TT_q(M_q) \simeq \TT_q(\La_{M_q})$.

Let $\B\in M_{N\times \ex}(\Zset)$ be an exchange matrix and $\La= (\la_{ij}) \in M_N(\Zset)$ be a skew-symmetric matrix.
The pair $(\La, \B)$ is called {\em{compatible}} if 
\[
\sum_{k=1}^N b_{kj}\lambda_{ki} = \delta_{ij} d_j, \quad \forall k \in [1,N], j \in \ex
\]
for a collection of 
%relatively prime 
positive integers $(d_j, j \in \ex)$. In terms of the diagonal matrix $D := \diag(d_j, j \in \ex)$, this condition is written as
\begin{equation}
\label{compatibility}
\B^\top \La = [D \; 0], 
\end{equation}
where where 0 denotes the zero matrix of size $\ex \times ([1, N] \backslash \ex)$.
If $(\La, \B)$ is a compatible pair, then $\B$ has full rank and its principal part $B$ is skew-symmetrized by $D$.

For each $k\in \ex$, the mutation of a compatible pair $(\La, \B)$ in the direction of $k$ is 
\begin{equation}
\label{mut-com-pairs}
\mu_k(\La, \B) := (\La',\B'),
\end{equation} 
where  $\B' = E_s \B F_s$ as in Sect. \ref{2.2} and $\La':= E_s^\top \La E_s$, $s = \pm$. 
The pair $(\La',\B')$ is compatible with respect to the same diagonal matrix $D$ 
and is independent on the choice of sign $s$.

A {\em{quantum seed}} is a pair $(M_q, \B)$, consisting of a toric frame $M_q$ of $\FF_q$ and an exchange matrix $\B$
such that $(\La_{M_q}, \B)$ is compatible.
We call $M_q(e_j)$, $j\in [1,N]$, \emph{cluster variables} of the seed, among which the {\em{frozen}} ones are those indexed by $[1,N]\backslash \ex$.
Denote $b^k$ be the $k$-th column of $\B$. The mutation of a quantum seed $(M_q, \B)$ in the direction of $k \in \ex$ is 
\[
\mu_k(M_q, \B) = (\mu_k(M_q), \mu_k(\B) ) := (\rho^{M_q}_{b^k, s} M_q E_s , E_s \B F_s ) 
\]
for any choice of sign $s$, 
where $\rho_{b^k, s} = \rho^{M_q}_{b^k, s}$ is the unique automorphism of $\FF_q$ such that 
\[
\rho_{b^k, s}(M_q(e_j)) = 
\begin{cases}
M_q(e_k) + M_q(e_k + s b^k) & \text{if } j=k \\
M_q(e_j)                   &      \text{if } j \neq k.
\end{cases}
\]
The skew-symmetric matrix associated to the toric frame $\mu_k(M_q)$ is $\La_{\mu_k(M_q)} = \mu_k(\La_{M_q})$.
The mutation process is involutive.  
 
The {\em{quantum cluster algebra}} $\A_q(M_q, \B, \inv)$ is the $\Abb$-subalgebra of $\FF_q$ generated by all cluster variables of quantum seeds $(M_q',\B')\sim (M_q, \B)$ 
and by the inverted frozen variables $\{ M_q(e_j)^{-1} \mid j \in \inv \}$.
The corresponding \emph{upper quantum cluster algebra} is given by the intersection
\[
\U_q(M_q, \B, \inv):= \bigcap_{(M_q',\B')\sim (M_q,\B)} \TT_q(M_q')_\geq
\] 
of the \emph{mixed quantum tori}
\[
\TT_q(M'_q)_{\geq}:=\Abb\lcor \hspace{1pt} M_q'(e_k)^{\pm 1}, M_q'(e_i) 
\mid k \in \ex \sqcup \inv, 
i \in \ninv \rcor 
\subset \TT_q(M'_q).
\]
Here, the term mixed refers to the fact that these algebras are mixtures of quantum tori and quantum affine spaces.
We have the quantum Laurent phenomenon:
\begin{equation}
\label{qLaurent}
\A_q(M_q, \B, \inv) \subseteq \U_q(M_q, \B, \inv),
\end{equation}
proved in \cite{BerZe} for $\inv = [1,N] \backslash \ex$ and \cite{GY1} in general.

The {\em{exchange graphs}} of the upper cluster algebra $\U(\B, \inv)$ and the upper quantum cluster algebra $\U_q(M_q, \B, \inv)$ 
are the labelled graphs with vertices corresponding to the seeds that are mutation-equivalent to $(\wt{\mathbf{x}}, \B)$ and $(M_q, \B)$, respectively,
and edges given by seed mutation, labelled by the corresponding mutation number. (The exchange graph does not depend on the set $\inv$ of inverted frozen variables.) 
\bth{isom-e-graphs} {\em{(Berenstein--Zelevinsky)}} \cite[Theorem 6.1]{BerZe}
There exists a unique isomorphism of labelled graphs, between the exchange graphs of $\U(\B, \inv)$ and $\U_q(M_q, \B, \inv)$ 
that sends the vertex corresponding to the seed $(\wt{\mathbf{x}}, \B)$ to that of the seed $(M_q, \B)$.
\eth 
%%%%%%%%%%%%
\subsection{Poisson structures on cluster algebras}
\label{2.4}
We follow Gekhtman, Shapiro and Vainshtein \cite{GSV}. Consider an upper cluster algebra $\U(\B, \inv)$ for which there exists a skew-symmetric matrix $\La \in M_N(\Zset)$ 
such that $(\La, \B)$ is a compatible pair. This is equivalent to saying that there exists a quantum cluster algebra with a seed with exchange matrix $\B$.  
For such a seed $(M_q, \B)$, the skew-symmetric matrix of the toric frame $M_q$ equals $\La$.

In this setting, for every seed $(\x',\B') \sim (\x, \B)$, there exists a skew symmetric matrix $\La_{(\x',\B')}$ satisfying the following two conditions:
\begin{enumerate}
\item[(i)] the pair $\big( \La_{(\x', \B')}, \B' \big)$ is compatible and
\item[(ii)] for all $k \in \ex$, 
\[
\mu_k \big( \La_{(\x', \B')}, \B' \big) = \big( \La_{(\mu_k (\x') , \mu_k(\B'))}, \mu_k(\B')\big)
\]
where the left hand side uses mutation of compatible pairs, see \eqref{mut-com-pairs}. 
\end{enumerate}
For a seed  $(\x',\B') \sim (\x, \B)$, we choose a sequence of mutations such that 
\[
\mu_{k_1} \ldots \mu_{k_j} (\x, \B) = (\x', \B').
\]
Applying the sequence of mutations to the compatible pair $(\La, \B)$ gives that
\[
\mu_{k_1} \ldots \mu_{k_j} (\La, \B) = (\La', \B')
\]
for a skew-symmetric matrix $\La' \in M_N(\Zset)$. The exchange graph isomorphism from \thref{isom-e-graphs} implies that the matrix $\La'$ is independent 
on the choice of mutation sequence $\mu_{k_1} \ldots \mu_{k_j}$. We set $\La_{(\x', \B')} := \La'$. 
It is clear that conditions (i)-(ii) are satisfied. 

For each seed $(\x', \B')\sim (\x, \B)$, the mixed polynomial/Laurent polynomial ring $\LL(\x')$ has a Poisson algebra structure such that
\begin{equation}
\label{Poisson-L}
\{ x'_k, x'_i \} = \La_{(\x', \B')}(e_k,e_i) x'_k x'_i, \quad \forall k,i \in [1,N].
\end{equation}
All algebras  $\LL(\x')$ have a common field of fractions $\FF$, so all Poisson brackets \eqref{Poisson-L} automatically extend to Poisson field structures on $\FF$.
Those extensions coincide and, hence, the intersection \eqref{upper-c-a},  $\U(\B, \inv)$, inherits a Poisson algebra structure, 
called the Gekhtman--Shapiro--Vainshtein (GSV) Poisson algebra structure of $\U(\B, \inv)$. 
%We will refer to it as to the GSV Poisson structure and  denote it by $\{.,.\}_{\GSV}$. 
%%%%%%%%%%%%%%%%%%%%%%%%%%%%%%%
\sectionnew{Root of unity quantum cluster algebras and torus automorphisms}
In this section we gather background material on root of unity quantum cluster algebras from \cite{NTY2,HLY} and prove auxiliary properties
that will be used in the next sections. In particular, we describe certain torus actions on these algebras, which will play a key role in the paper.
%%%%%%%%%%%%%%%%%%%%%%
\subsection{Root of unity quantum cluster algebras}
\label{3.1}
Let $\ell$ be a positive integer and $\ep^{1/2}\in \Cset$ be a primitive $\ell$-th root of unity. Denote
\[
\Zset/\ell:= \Zset/{\ell\Zset}. 
\]
In this paper we will work over $\Cset$, while the construction in \cite{NTY2} was carried out over $\Zset[\ep^{1/2}]$.
We start with a skew-symmetric bilinear form $\Omega : \Zset^N \times \Zset^N \to \Zset/\ell$
and identify it with the skew-symmetric matrix $( \Omega(e_i, e_j))_{i,j=1}^N \in M_N(\Zset/\ell)$. 
The associated {\em root of unity based quantum torus} is
\begin{equation}
\label{Tep}
\TT_\ep(\Omega) := \Span_{\Cset} \{ x^f \mid f \in \Zset^N \}, \; \; 
\mbox{where} \; \; x^f x^g = \ep^{\Omega(f,g)/2} x^{f+g}, \; \; \forall f,g \in \Zset^N.
\end{equation}

A {\em root of unity toric frame} $M_\ep$ of a division algebra $\FF_\ep$ over $\Cset$ is a map 
\[
M_\ep : \Zset^N \to \FF_\ep
\]
for which there exists an $\Cset$-algebra embedding $\phi: \TT_\ep(\Omega) \hra \FF_\ep$ for a skew-symmetric matrix $\Omega\in M_N(\Zset/\ell)$ 
with the properties that $\phi(x^f)=M_\ep(f)$ for all $f\in \Zset^N$ and $\FF_\ep \simeq \Fract \left( \TT_\ep(\Omega) \right)$.
Denote by $\Omega_{M_\ep}$ the skew-symmetric matrix (bilinear form) with values in $\Zset/\ell$ of the root of unity toric frame $M_\ep$, 
given by
\[
M_\ep(f) M_\ep(g) = \ep^{\Omega(f,g)/2} M_\ep(f+g), \quad \forall f, g \in \Zset^N.
\] 
The mixed quantum tori for the root of unity setting are given by
\begin{equation}
\label{Tep-geq}
\TT_\ep(M_\ep)_{\geq} :=\Cset \lcor M_\ep(e_k)^{\pm 1}, \ M_\ep(e_i) \mid k \in \ex \sqcup \inv,  i \in \ninv \rcor \subset \TT_\ep( M_\ep).
\end{equation}

A pair $(\Omega, \B)$, consisting of a skew-symmetric matrix $\Omega  \in M_N(\Zset/\ell)$ and an exchange matrix $\B \in M_{N \times \ex}(\Zset)$,
is said to be {\em{$\ell$-compatible}} if there exists a diagonal matrix $D := \diag(d_j, j \in \ex)$ with $d_j \in \Zset_+$ 
such that the principal part $B$ of $\B$ is skew-symmetrized by $D$ and 
\[
\ol{\B}^\top \hspace{-0.15cm} \Omega  = 
\begin{bmatrix}
\, \ol{D} \; 0
\end{bmatrix} \hspace{-0.1cm}.
\]
Here and below, for an integer matrix $C$, $\ol{C}$ denotes its reduction modulo $\ell$.
%Since there is no restriction on $\ol{d}_j$, so the matrix $\B$ need not have full rank as seen in the case of quantum cluster algebras. 

The \emph{mutation} in direction $k \in \ex$ of an $\ell$-compatible pair is defined to be $\mu_k ( \Omega, \B) := ( \ol{E}_s^\top \Omega \ol{E}_s, E_s \B F_s)$ for $s= \pm$, 
and as in the quantum case; it is independent of the choice of sign $s$. Further, the pair
$\mu_k(\Omega, \B)$ is also $\ell$-compatible with respect to the same diagonal matrix $D$. %The mutation of $\ell$-compatible pairs is an involution.

A pair $(M_\ep, \B)$ is called a {\em root of unity quantum seed} if $(\Omega_{M_\ep},\B)$ is an $\ell$-compatible pair. The mutation in direction $k\in\ex$ of a root of unity quantum seed is similar to that in the quantum seed case,
$\mu_k(M_\ep, \B) := (\rho^{ M_\ep}_{b^k, s} M_\ep E_s,  E_s \B F_s)$, where $\rho_{b^k,s}=\rho^{M_\ep}_{b^k,s}$ is the unique automorphism of $\FF_\ep$ given by 
\[ \rho_{b^k, s}^{M_\ep}(M_\ep(e_j)) = 
\begin{cases}
M_\ep(e_k) + M_\ep(e_k + sb^k) & \text{if } j=k \\
M_\ep(e_j)                & \text{if } j \neq k.
\end{cases}\]
The skew-symmetric matrix associated to $\mu_k(M_\ep)$ is $\Omega_{\mu_k(M_\ep)} = \mu_k(\Omega_{M_\ep})$. 
Moreover, mutation of root of unity quantum seeds does not depend on the choice of sign and is an involution. 
Two seeds are called {\em{mutation equivalent}}, $(M'_\ep, \B')\sim (M''_\ep, \B'')$, if one can be obtained from the other by a finite sequence of mutations.

The {\em root of unity quantum cluster algebra} $\A_\ep(M_\ep, \B, \inv)$ is defined to be the $\Cset$-subalgebra of $\FF_\ep$ 
generated by all cluster variables of the seeds $(M'_\ep, \B') \sim (M_\ep, \B)$ and by 
$M_\ep(e_i)^{-1}$ for $i \in \inv$.
The corresponding \emph{root of unity upper quantum cluster algebra} 
is defined as the intersection 
\[
\U_\ep(M_\ep, \B, \inv) := \bigcap_{ (M'_\ep, \B') \sim (M_\ep, \B)}\TT_\ep(M_\ep')_{\geq}.
\]
We have $\A_\ep(M_\ep, \B, \inv)\subseteq \U_\ep(M_\ep,  \B, \inv)$ by \cite[Theorem 3.10]{NTY2}.

For every root of unity toric frame $M'_\ep$ of $\FF_\ep$ and $1 \leq i \leq N$, 
\[
M'_\ep(e_i)^\ell \in \ZZ(\FF_\ep). 
\] 
Consider the central subalgebra
\begin{equation}
\label{Lep}
\LL_\ep(\Omega') = \Cset [ (x^{e_k})^{\pm \ell}; 1 \leq k \leq N] \subset \ZZ(\TT_\ep(\Omega'))
\end{equation}
and the mixed polynomial/Laurent polynomial central subalgebras
\begin{align}
\label{Lep-geq0}
\LL_\ep( \Om')_{\geq} &:= \Cset [ (x^{e_k})^{\pm \ell}, \ (x^{e_i})^\ell  ;  
k \in \ex \sqcup \inv,  i \in \ninv ] \subset \ZZ( \TT_\ep(\Om')_\geq), \\
\label{Lep-geq}
\LL_\ep( M'_\ep)_{\geq} &:= \Cset [ M'_\ep(e_k)^{\pm \ell}, \ M'_\ep(e_i)^\ell  ;  
k \in \ex \sqcup \inv,  i \in \ninv ] \subset \ZZ( \TT_\ep(M'_\ep)_\geq).
\end{align}
Denote the central subalgebra
\begin{equation}
\label{CU}
\C\U_\ep(M_\ep, \B, \inv) := \bigcap_{ (M'_\ep, \B') \sim (M_\ep, \B)} \LL_\ep(M_\ep')_{\geq}
\subset \ZZ (\U_\ep(M_\ep, \B, \inv)).
\end{equation}
The {\em{exchange graph}} of the root of unity upper quantum cluster algebra $\U_\ep(M_\ep, \B, \inv)$ is 
the labelled graph with vertices corresponding to the root of unity quantum seeds that are mutation-equivalent to $(M_\ep, \B)$
and edges given by seed mutations, labelled by the corresponding mutation number. (The exchange graph is independent of the set $\inv$.)

\bth{cent-isom} Assume that $\ell$ is an odd positive integer that is coprime to the diagonal entries of the skew-symmetrizing matrix $D$ 
for the principal part of the exchange matrix $\B$. Then the following hold:
\begin{enumerate}
\item[(i)] \cite[Theorem 4.8]{NTY2}
The exists a unique isomorphism of labelled graphs, between the exchange graphs of $\U(\B, \inv)$ and $\U_\ep(M_\ep, \B, \inv)$ 
that sends the vertex corresponding to the seed $(\wt{\mathbf{x}}, \B)$ to that of the seed $(M_\ep, \B)$.
\item[(ii)] \cite[Proposition 4.4]{NTY2} For all seeds $(M_\ep', \B') \sim (M_\ep, \B)$ and $k \in \ex$,
\[
M'_\ep(e_k)^\ell \left( \mu_k M'_\ep (e_k) \right)^{\ell} = \prod_{b'_{ik}>0} (M'_\ep(e_i)^\ell)^{b'_{ik}} + \prod_{b'_{ik}<0} (M'_\ep(e_i)^\ell)^{ - b'_{ik}},
\]
which is precisely the mutation formula for the seeds of $\U(\B, \inv)$.
\item[(iii)] \cite[Proposition 3.9]{HLY} We have an isomorphism of $\Cset$-algebras 
\begin{equation}
\label{U-CU-isom}
\C\U_\ep(M_\ep, \B, \inv) \simeq \U(\B, \inv),
\end{equation}
which is uniquely determined by sending $M'_\ep(e_k)^\ell \mt x'_k$, where the seed $(M'_\ep, \B')$ corresponds 
to the seed $(\x', \B')$ under the graph isomorphism in part (i).  
\item[(iv)] \cite[Theorem B]{HLY} $\U_\ep(M_\ep, \B, \inv)$ is a Cayley--Hamilton algebra of degree $\ell^N$ over its central subalgebra $\C\U_\ep(M_\ep, \B, \inv)$
in the sense of Procesi \cite[Definition 4.2]{DP}. 
$\U_\ep(M_\ep, \B, \inv)$ is finitely generated $\Cset$-algebra if and only if $\U_\ep(M_\ep, \B, \inv)$ is a finitely generated 
$\C\U_\ep(M_\ep, \B, \inv)$-module and $\U(\B, \inv)$ is a finitely generated $\Cset$-algebra. 
\end{enumerate}
\eth
We will need the root of unity upper quantum cluster algebra in the special case when all frozen variables are inverted:
\[
\U_\ep(M_\ep, \B):= \U_\ep(M_\ep, \B, [1,N]\backslash \ex) = \bigcap_{(\x', \B') \sim (\x, \B)} \TT_\ep(M_\ep').
\]
It is easy to verify that it is a localization of the ones for other choices of $\inv$:
\[
\U_\ep(M_\ep, \B) = \U_\ep(M_\ep, \B, \inv)[M(e_i)^{-1}; i \in \ninv].
\]

The above treatment can be viewed as defining quantum cluster $\mathcal{A}$-varieties at roots of unity. 
Quantum cluster $\mathcal{X}$-varieties at roots of unity were defined and studied by Fock and Goncharov in \cite{FG}, who obtain analogous algebraic results
to \thref{cent-isom}(ii)-(iii) under the stronger assumption that  the order of the root of unity is coprime to the entries of the exchange matrices of all seeds of the algebra.
However, they only consider varieties up to birational isomorphism and therefore do not consider such phenomena as singularities.
\subsection{Strict root of unity quantum cluster algebras}
\label{3.2}
\bde{strict-seed} \cite[Sect. 5]{NTY2} We say that a root of unity quantum seed $(M_\ep, \B)$ is {\em{strict}} if there exists a skew-symmetric integer matrix 
$\La \in M_N (\Zset)$ such that 
\begin{enumerate} 
\item[(i)] $\Omega_{M_\ep} = \overline{\La}$ and 
\item[(ii)] $(\La, \B)$ is a compatible pair, 
%with respect to the diagonal matrix $D$, which skew-symmetrizes the principal part of $\B$ 
see Sect. \ref{2.3}.
\end{enumerate}
The corresponding root of unity upper quantum cluster algebra $\U_\ep(M_\ep, \B, \inv)$ will be also called strict. 
\ede
\bpr{struct-mut} Assume that $\U_\ep(M_\ep, \B, \inv)$ is a strict root of unity upper quantum cluster algebra such that $\ell$ 
is an odd positive integer that is coprime to the diagonal entries of the skew-symmetrizing matrix $D$ for the principal part of the exchange matrix $\B$. 
Then for every seed $(M'_\ep, \B') \sim (M_\ep, \B)$ there exists a unique skew-symmetric integer matrix 
$\La_{(M'_\ep, \B')} \in M_N (\Zset)$, such that $\La_{(M'_\ep, \B')} = \La'$, the matrix from \deref{strict-seed}, and 
\begin{enumerate}
\item[(i)] the pair $\big( \La_{(M'_\ep, \B')}, \B' \big)$ is compatible and $\overline{\La_{(M'_\ep, \B')}} = \Omega_{M'_\ep}$;
%pair with respect to the diagonal matrix $D$ skew-symmetrizing the principal part of $\B$ and
\item[(ii)] for all $k \in \ex$, 
\[
\mu_k \big( \La_{(M'_\ep, \B')}, \B' \big) = \big( \La_{(\mu_k (M'_\ep, \B') )}, \mu_k(\B') \big),
\]
where the left hand side uses mutation of compatible pairs, defined in \eqref{mut-com-pairs}. 
\end{enumerate}
\epr
For the proof of the proposition we will need the surjective $\Cset$-algebra homomorphism
\begin{equation}
\label{kappa}
\kappa_\ep : \TT_q(M_q)_\geq \to \TT_q(M_\ep)_\geq,
\end{equation}
given by
\[
q^{1/2} \mt \ep^{1/2}, M_q(e_k)^{\pm 1} \mt M_\ep(e_k)^{\pm 1},
M_q(e_i) \mt M_\ep(e_i), \forall k \in \ex \sqcup \inv, i \in \ninv.
\]
Its kernel is 
\begin{equation}
\label{ker-kappa}
\ker \kappa_\ep = (q^{1/2} - \ep^{1/2}) \TT_q(M_q)_\geq,
\end{equation}
see \cite[Lemma 5.5]{NTY2}.  
\medskip

\noindent
{\em{Proof of \prref{struct-mut}.}}
For a given seed  $(M'_\ep, \B') \sim (M_\ep, \B)$ consider a sequence of mutations $\mu_{k_1} \ldots \mu_{k_j}$ such that 
\begin{equation}
\label{mut-sequence}
\mu_{k_1} \ldots \mu_{k_j} (M_\ep, \B) = (M'_\ep, \B').
\end{equation}
Applying this sequence of mutations to the compatible pair $(\La, \B)$ gives that
\begin{equation}
\label{comp-pair}
\mu_{k_1} \ldots \mu_{k_j} (\La, \B) = (\La', \B')
\end{equation} 
for a skew-symmetric matrix $\La' \in M_N(\Zset)$. The two exchange graph isomorphisms from Theorems 
\ref{tisom-e-graphs} and \ref{tcent-isom}(i) imply that the matrix $\La'$ is independent on the choice of mutation sequence $\mu_{k_1} \ldots \mu_{k_j}$
satisfying \eqref{mut-sequence}. Define
\[
\La_{(M'_\ep, \B')} := \La'.
\]
Since the mutations of a compatible pair are compatible pairs (see Sect. \ref{2.3}), $\big( \La_{(M'_\ep, \B')}, \B' \big)$
is a compatible pair. 
By \cite[Theorem 5.7]{NTY2}, $\kappa_\ep$ restricts to a surjective $\Cset$-algebra homomorphism $\A_q(M_q, \B, \inv) \to \A_\ep(M_\ep, \B, \inv)$ 
that sends cluster variables to cluster variables and commutes with mutation. The toric frame 
$M_q$ of $\FF_q$ is such that $\La_{M_q} = \La$. The existence of such a homomorphism implies that
\[
M'_\ep(e_i) M'_\ep(e_k) = \ep^{\La'(e_i, e_k)/2} M'_\ep(e_i+e_k), \quad \forall i,k \in [1,N].
\]
%for every seed $(M'_\ep, \B)$ of $\U_\ep(M_\ep, \B, \inv)$.
Hence, 
\[
\overline{\La_{(M'_\ep, \B')}} = \Omega_{M'_\ep}.
\]
Finally, property (ii) in the statement of the proposition follows at once from the definition \eqref{comp-pair} of $\La_{(M'_\ep, \B')}$.
\qed
%%%%%%%%%%%%%%%%%%%%%%
\subsection{Torus action on $\U_\ep(M_\ep, \B, \inv)$}
\label{3.3}
Denote by $\Ker(\B^\top)$ the null space of $\B^\top$ in $\Zset^N$. 
Let $\nu=(\nu_1, \dots, \nu_N) \in \Ker(\B^\top)$. For each $k\in \ex$, 
we define the vector $\mu_k(\nu)=(\nu_1, \dots,\nu_{k-1}, \nu_k', \nu_{k+1}, \dots, \nu_N)$, where
 \[\nu_k' := \nu\cdot[b^k]_+ - \nu_k.\]
Recall that $b^k$ denotes the $k$-th column of $\B$ and $f \cdot g$ the dot product on $\Zset^N$. 

By \cite[Lemma 2.3]{GSV},
\begin{equation}
\label{ker-nu}
\nu\in \Ker(\B^\top) \quad \Rightarrow \quad 
\mu_k(\nu)\in \Ker(\mu_k(\B)^\top), \; \; \forall k \in \ex.
\end{equation}
For $\nu \in  \Ker(\B^\top)$, we have the $\Cset^\times$-action on the root of unity mixed quantum torus $\TT_\ep(M_\ep, \B)_\geq$, given by 
\begin{equation}
\label{phi-act}
\varphi_\nu(t) \cdot M_\ep(f):= t^{\nu \cdot f} M_\ep(f), \; \; \forall f \in \Zset^N.
\end{equation}
It induces $\Cset^\times$-action on the skew field of fractions $\FF_\ep$ of $\TT_\ep(M_\ep, \B)$. 
Analogously to the proof of \cite[Lemma 2.3]{GSV} one verifies that
\begin{equation}
\label{phi-nu-mut}
\varphi_\nu(t)\cdot \mu_k(M_\ep)(f) := t^{\mu_k(\nu) \cdot f} \mu_k(M_\ep)(f), \; \; \forall f \in \Zset^N.
\end{equation}
By recursively applying \eqref{ker-nu}, we obtain that this action preserves the root of unity mixed quantum tori $\TT_\ep(M'_\ep, \B')_\geq$
associated to all seeds $(M'_\ep, \B') \sim (M_\ep, \B)$, and thus, the action preserves $\U_\ep(M_\ep, \B, \inv)$.
This proves the first part of the following lemma:

\ble{Caction} 
\hfill
\begin{enumerate}
\item[(i)]
For every $\nu \in \Ker(\B^\top)$, there is a $\Cset^\times$-action $\varphi_\nu$ on the root of unity upper quantum 
cluster algebra $\U_\ep(M_\ep, \B, \inv)$, given by \eqref{phi-act}. For every seed $(M'_\ep, \B') \sim (M_\ep, \B)$ and $k \in [1,N]$, 
\begin{equation}
\label{phi-nu-M'}
\varphi_\nu(t) \cdot M'_\ep(e_k) = t^a M'_\ep(e_k)
\end{equation}
for some $a \in \Zset$ depending on the seed and $k$. 
\item[(ii)] The action $\varphi_\nu$ preserves  $\C\U_\ep(M_\ep, \B, \inv)$. Under the isomorphism \eqref{U-CU-isom}, 
it corresponds to the $\Cset^\times$-action $\psi_\nu$ on $\U(\B, \inv)$ uniquely determined from 
\[
\psi_\nu(t) \cdot x_k := t^{\ell (\nu \cdot e_k) } x_k, \; \; \forall k \in [1,N].
\]
For every seed $(\x', \B') \sim (\x, \B)$ and $k \in [1,N]$, 
\[
\psi_\nu(t) \cdot x'_k := t^{\ell a } x'_k 
\]
for some $a \in \Zset$.
\item[(iii)] For every $\nu, \eta \in \Ker(\B^\top)$, the $\Cset$-actions $\varphi_\nu$ and $\varphi_\eta$ commute.
\end{enumerate}
\ele
The second and third parts of the lemma follow from \eqref{phi-nu-M'} and the concrete form of the isomorphism \eqref{U-CU-isom}
from \thref{cent-isom}(iii). Up to the $\ell$-th power the $\Cset$-action $\psi_\nu$ on  $\U(\B, \inv)$ is the one constructed in 
\cite[Lemma 2.3]{GSV}. 

\bde{nB} For an exchange matrix $\B$, define the nullity
\[
n(\B) := \dim \Ker(\B^\top)
\]
and the complex torus
\[
T(\B) := (\Cset^\times)^{n(\B)}.
\]
\ede
$\Ker(\B^\top)$ is a free abelian group of rank $n(\B)$. Fix a basis 
\begin{equation}
\label{basis-KerB}
\{\nu^1, \dots, \nu^{n(\B)} \}
\end{equation}
of it. Parts (i) and (iii) \leref{Caction} imply that we have a 
$T(\B)$-action on $\U_\ep(M_\ep, \B, \inv)$, given by 
\[
\varphi(t_1, \dots, t_{n(\B)}) \cdot y := \varphi_{\nu^1}(t_1)\cdots \varphi_{\nu^{n(\B)}}(t_{n(\B)}) \cdot y, 
\]
for $t_i \in \Cset$ and $y \in \U_\ep(M_\ep, \B, \inv)$.
For every seed $(M'_\ep, \B') \sim (M_\ep, \B)$ and $k \in [1,N]$, there exists a character $\theta : T(\B) \to \Cset^\times$ such that 
\begin{equation}
\label{phi-M'}
\varphi(t) \cdot M'_\ep(e_k) = \theta(t) M'_\ep(e_k), \; \; \forall t \in T(\B).
\end{equation}
For different choices of a basis $\{\nu^1, \dots, \nu^m\}$ of $\Ker(\B^\top)$, the actions $\varphi$ differ from each other by an automorphism of $T(\B)$. 

Part (ii) of \leref{Caction} implies that the action $\varphi$ preserves $\C\U_\ep(M_\ep, \B, \inv)$. Under the isomorphism \eqref{U-CU-isom}, 
it corresponds to the $T(\B)$-action $\psi$ on $\U(\B, \inv)$ given by 
\[
\psi(t_1, \dots, t_{n(\B)}) \cdot y := \psi_{\nu^1}(t_1)\cdots \psi_{\nu^{n(\B)}}(t_{n(\B)}) \cdot y, 
\]
for $t_i \in \Cset$ and $y \in \U(\B, \inv)$.
For every seed $(\x', \B') \sim (\x, \B)$ and $k \in [1,N]$, there exists a character $\theta : T(\B) \to \Cset^\times$ such that 
\begin{equation}
\label{rescale-x}
\psi(t) \cdot x'_k = \theta(t^\ell) x'_k, \; \; \forall t \in T(\B).
\end{equation}
By the results in \cite{GSV}, the GSV Poisson algebra structure on $\U(\B, \inv)$ is invariant under the 
action $\psi$ of $T(\B)$. The action $\psi$ is defined for every upper cluster algebra $\U(\B, \inv)$.

Finally, we note that the $T(\B)$-action $\psi$ on $\U(\B, \inv)$ is defined for any upper cluster algebra, without the need of 
existence of root of unity quantization $\U_\ep(M_\ep, \B, \inv)$
%%%%%%%%%%%%%%%%%%%
\sectionnew{Poisson side: the open torus orbit of symplectic leaves}
\lb{Torbits} 
In this section we prove our first main result proving that the GSV Poisson structure on every finitely generated 
upper cluster algebra always has a Zariski open torus orbit of symplectic leaves and explicitly describe this set.  
\subsection{Statement of main theorem}
\label{5.1}
Throughout this section, we will assume that $\U(\B, \inv)$ is a finitely generated, upper cluster algebra for which there exists a skew-symmetric matrix $\La \in M_N(\Zset)$ 
such that $(\La, \B)$ is a compatible pair. The GSV Poisson algebra structure $\{.,.\}$ on $\U(\B, \inv)$ (see Sect. \ref{2.4})
gives rise to a Poisson structure $\pi$ on the affine variety 
\[
Y(\B):= \MaxSpec \U(\B, \inv).
\]
This variety is in general singular. It is normal because each mixed polynomial/Laurent polynomial ring $\LL(\x')_\geq$ is integrally closed, 
and thus the upper cluster algebra $\U(\B, \inv)$, given by the intersection \eqref{upper-c-a}, is integrally closed too. 

Denote the product of non-inverted frozen variables
\begin{equation}
\label{x-prod}
x := \prod_{i \in \ninv} x_i \in \Cset[Y(\B)].
\end{equation}
Denote
\[
Y(\B)^\circ = Y(\B) \backslash \VV(x) \cong \MaxSpec \U(\B).
\]
For a seed $(\x', \B') \sim (\x, \B)$, denote
\begin{align*}
S_{\x'} &:=\MaxSpec \LL(\x')_{\geq}  \cong (\Cset^\times)^{ | \ex | + | \inv |} \times \Cset^{| \ninv |}, \\ 
S_{\x'}^\circ &:=\MaxSpec \LL(\x') \cong (\Cset^\times)^N.
\end{align*}
Clearly, 
\[
S_{\x'}^\circ \cong S_{\x'} \backslash \VV(x).
\]
The algebras $\LL(\x')_{\geq}$ and $\LL(\x')$ are localizations of $\U(\B, \inv)$ by the multiplicative subsets generated by $x_i$ for $i \in \ex$
and $i \in \ex \sqcup \ninv$, respectively. So,  $S_{\x'}$ and $S_{\x'}^\circ$ are Zariski open, and thus, dense subsets of $Y(\B)$. 

Recall from Sect. \ref{3.3} that the Poisson structure $\pi$ is invariant under the action $\psi$ of $T(\B)$. Thus it makes sense to consider the 
$T(\B)$-orbits of symplectic leaves of $(\U(\B, \inv), \pi)$, which are regular Poisson submanifolds.
\bth{open-T-orb-sl} Assume that $\U(\B, \inv)$ is a finitely generated upper cluster algebra for which there exists a skew-symmetric matrix $\La \in M_N(\Zset)$ 
such that $(\La, \B)$ is a compatible pair. Then the affine Poisson variety $(Y(\B), \pi)$ has a Zariski open $T(\B)$-orbit of symplectic leaves, which equals
\[
Y(\B)^\reg \backslash \VV(x) = (Y(\B)^\circ)^\reg.
\]
\eth
%%%%%%%%%%%%%%%%%%%%
\subsection{One-sided containment}
\label{5.2}
\bde{P-prime} \cite[Sect. 4]{GY3}, \cite[Sect. 2.2]{NTY1} Let $R$ be a Poisson algebra which is an integral domain considered as a commutative algebra. 
An element $p \in R$ is called Poisson prime, if the principal ideal $(p)$ is a Poisson ideal and a prime ideal of $R$, considered 
as a commutative algebra.

Equivalently, $p \in R$ is a prime element and 
\[
p \mid \{ p, r \}, \quad \forall r \in R. 
\]
\ede
\bpr{p-primes} Assume that $\U(\B, \inv)$ is an upper cluster algebra, such that $\B$ admits a compatible skew-symmetric matrix $\La \in M_N(\Zset)$. 
Then the following hold:
\begin{enumerate}
\item[(i)] Every prime element of $\Cset[x_i; i \in \ninv]$ is a prime element of $U(\B, \inv)$.
\item[(ii)] Every non-inverted frozen variable $x_i$, {\em{(}}i.e. $i \in \ninv${\em{)}}
is a Poisson prime of $\U(\B, \inv)$ with respect to the Poisson structure associated to $\La$. 
\end{enumerate} 
\epr
\begin{proof} (i) Let $p \in \Cset[x_i; i \in \ninv]$ be a prime element. Analogously to the proof of \cite[Proposition 3.5]{GLS-spec}, one shows that
\begin{equation}
\label{inters}
\LL(\x')_{\geq} \cap ( p \LL(\mu_k \x')_{\geq}) = ( p \LL(\x')_{\geq}) \cap \LL(\mu_k \x')_{\geq}
\end{equation}
for all seed $(\x', \B') \sim (\x,\B)$ and $k \in \ex$. Assume that $p | ab$ for some $a, b \in U(\B, \inv)$. Then $a, b \in \LL(\x)_\geq$, 
and since $p$ is a prime element of $\LL(\x)_\geq$, either $p \mid a$ or $p \mid b$. Say $p \mid a$. 
Then $a \in p \LL(\x)_{\geq}$, and by \eqref{inters}, 
\[
a \in p \LL(\mu_k \x)_{\geq}, \quad \forall k \in \ex.
\]
Iterating this argument gives that $a \in p \LL(\x')_{\geq}$ for all $(\x', \B') \sim (\x,\B)$. Hence,
\[
a \in \bigcap_{(\x', \B') \sim (\x,\B)}  p \LL(\x')_{\geq} = p \, U(\B, \inv).
\] 
This implies the statement of (i) because $p$ is not a unit of $U(\B, \inv)$.

(ii) Fix $b \in U(\B, \inv)$. It follows from \eqref{Poisson-L} that
\[
\{x_i, b \} /x_i \in \LL(\x') 
\]
for all $(\x', \B') \sim (\x,\B)$. Hence, 
\[
\{x_i, b \} /x_i \in U(\B, \inv), 
\]
which completes the proof of part (ii). 
\end{proof}
\bpr{containment} In the setting of \thref{open-T-orb-sl}, $Y(\B)^\reg \backslash \VV(x) = (Y(\B)^\circ)^\reg$ is a union of $T(\B)$-orbits of symplectic leaves.
\epr
\begin{proof} Let $i \in \ninv$. By \eqref{rescale-x}, $x_i$ is rescaled by the $T(\B)$-action. So $\VV(x_i)$ is $T(\B)$-stable.  
It follows from  \cite[Proposition 2.3]{GY3} that $\VV(x_i)$ is a union of symplectic 
leaves. Therefore, $\VV(x_i)$ is a union of $T(\B)$-orbits of symplectic leaves. 

Since $T(\B)$ acts on $U(\B, \inv)$ by algebra automorphisms, $Y(\B)^\sing$ is $T(\B)$-stable. By \cite[Corollary 2.4]{P}, it is also 
a union of symplectic leaves. Thus, $Y(\B)^\sing$ is a union of $T(\B)$-orbits of symplectic leaves. 

Hence, $Y(\B)^\sing \cup \VV(x)$ is a union of $T(\B)$-orbits of symplectic leaves, and the same applies to its complement,
\[
Y(\B)^\reg \backslash \VV(x) = Y(\B) \backslash (Y(\B)^\sing \cup \VV(x)).
\qedhere
\]
\end{proof} 
%%%%%%%%%%%%%%%%%%%%
\subsection{A section of the anticanonical bundle of $Y(\B)^\reg$}
\label{5.3}
For $\nu \in \Ker( B^\top)$, denote by
\[
w_\nu \; \; \mbox{the vector field on $Y(\B)$},
\]
which represents the infinitesimal action of the $\Cset^\times$-action $\psi_\nu$ on $Y(\B)$, 
cf. \leref{Caction}(ii).

The rank of the skew-symmetric matrix $\La = (\la_{ik}) \in M_N(\Zset)$ equals 
\[
\rk (\La) = 2 r
\] 
for some $r \in \Zset_+$. It follows from \eqref{compatibility} that
\begin{align}
\label{nLa-ineq}
2r + n(\B) &= \dim \Im (\La) + \dim \Ker (\B^\top) \geq \dim \Im (\B^\top \La) + \dim \Ker (\B^\top) 
\\
&= \dim \Im (\B^\top) + \dim \Ker (\B^\top) =N.
\nn
\end{align}
The Poisson structure $\pi$ on $Y(\B)$ has rank $r$ (meaning $\wedge^r \pi \neq 0$ and $\wedge^{r+1} \pi =0$).
This follows from the fact that 
its rank equals the rank of its restriction to the Zariski open subset $S_\x$ on $Y(\B)$, on which the Poisson structure is 
explicitly given by 
\[
\pi|_{S_\x} = \sum_{1 \leq i, k \leq N} \la_{ik} \frac{\partial}{\partial x_i} \wedge \frac{\partial}{\partial x_k} \cdot
\]  
Recall from \eqref{basis-KerB} that the action $\psi$ of the torus $T(\B)$ on $Y(\B)$ is defined from a basis $\{\nu^1, \ldots, \nu^{n(\B)}\}$ 
of $\Ker( B^\top)$.
Denote by $\Delta$ the set of $(N-2r)$-element subsets of $[1,n(\B)]$, recall \eqref{nLa-ineq}. 
For $\theta := \{i_1, \ldots, i_{N- 2r}\} \in \Delta$, 
denote the section of the anticanonical bundle of $Y(\B)$,
\[
\chi_\theta := w_{v^{i_1}} \wedge \ldots \wedge w_{\nu^{i_{N-2r}}} \wedge (\wedge^r \pi) \in \Ga(Y(\B), K^*_{Y(\B)}).  
\]
\bth{sections-antican} Assume that $\U(\B, \inv)$ is a finitely generated upper cluster algebra, such that $\B$ admits a compatible skew-symmetric matrix $\La \in M_N(\Zset)$. 
\begin{enumerate}
\item[(i)] There exist a unique up to rescaling non-zero global section $\chi$ of the anticanonical bundle $K_{Y(\B)}^*$ such that for every $(\x', \B') \sim (\x, \B)$,
\[
\chi|_{S_{\x'}} = c_{(\x', \B')} x'_1\ldots x'_N  \frac{\partial}{\partial x'_1} \wedge \cdots \wedge \frac{\partial}{\partial x'_N}
\] 
for some $c_{(\x', \B')} \in \Cset^\times$. 
\item[(ii)] For each $\theta \in \Delta$, either $\chi_\theta=0$ or $\chi_\theta$ is a non-zero scalar multiple of $\chi$. 
The latter is the case for at least one $\theta \in \Delta$. 
\end{enumerate}
\eth
\begin{proof} First we show that for all $\theta \in \Delta$ and $(\x', \B') \sim (\x, \B)$,
\begin{equation}
\label{chi-theta}
\chi_\theta|_{S_{\x'}} = c_{\theta, (\x', \B')}  x'_1\ldots x'_N  \frac{\partial}{\partial x'_1} \wedge \cdots \wedge \frac{\partial}{\partial x'_N}
\end{equation}
for some $c_{\theta, (\x', \B')}  \in \Cset$. \leref{Caction}(ii) implies that 
\[
\psi_\nu(t) \cdot x'_k := t^{\ell a_k } x'_k, \; \; \forall k \in [1,N] 
\]
for some $a_k \in \Zset$. Therefore, 
\[
w_\nu|_{S_{\x'}^\circ} = 
\sum_{k=1}^N \ell a_k x'_k \frac{\partial}{\partial x'_k} \cdot 
\]
This equality and the fact that
\[
\pi|_{S_{\x'}} = \sum_{1 \leq i, k \leq N} \La_{(\x', \B')}(e_i,e_k) x'_i x'_k \frac{\partial}{\partial x'_i} \wedge \frac{\partial}{\partial x'_k}
\]
imply \eqref{chi-theta}. Since $S_{\x'}$ is a Zariski open subset of the irreducible affine variety $Y(\B)$, the following are equivalent 
for $\theta \in \Delta$:
\begin{enumerate}
\item $\chi_\theta \neq 0$, 
\item $c_{\theta, (\x, \B)} \neq 0$,
\item $c_{\theta, (\x', \B')} \neq 0$ for all $(\x', \B') \sim (\x, \B)$.
\end{enumerate}

The proof of \cite[Lemma 2.4]{GSV} gives that $T(\B)$ acts transitively on the symplectic leaves of $(S_{\x'}^\circ, \pi|_{V})$. Furthermore, 
for $y \in Y(\B)^\reg$ the tangent space at $y$ of the orbit $T(\B) \cdot y$ is
\begin{equation}
\label{span-T}
T_y \big( T(\B) \cdot y \big) = \Span_\Cset \{ w_{\nu^1, y}, \ldots, w_{\nu^{N- 2r}, y} \}.
\end{equation}
Therefore, there exists $\theta_0 \in \Delta$ such that 
\[
c_{\theta_0, (\x, \B)} \neq 0,
\]
which, in view of the above equivalence, implies that 
\[
\chi := \chi_{\theta_0} 
\]
satisfies the conditions in the proposition. Furthermore, for all $\theta \in \Delta$, $\chi_\theta \chi_{\theta_0}^{-1}$
is a meromorphic function on $Y(\B)$, which by \eqref{chi-theta} is constant on every Zariski open subset $S_{\x'}$, and thus, 
it should be constant on  $Y(\B)$. The same argument proves the uniqueness of a section $\chi$ with stated properties in part (i). 
\end{proof}
%%%%%%%%%%%%%%%%%%%%
\subsection{Normality and codimension two results}
\label{5.4}
Recall that an upper cluster algebra $\U(\B)$ is called totally coprime \cite[Sect. 1.2]{BFZ05} if every two columns of all of its extended exchange matrices are linearly independent. 
If the exchange matrix $\B$ has full rank, then this condition is satisfied. In turn, this condition is 
satisfied whenever $\B$ has a compatible skew-symmetric matrix $\La \in M_N(\Zset)$, \cite[Proposition 1.8]{BFZ05}.

\bpr{codim2} If $\U(\B)$ is a finitely generated, totally coprime upper cluster algebra, then 
\[
(Y(\B)^\circ)^\sing \quad \mbox{and} \quad 
Y(\B)^\circ \backslash \big( S_{\x} \cup (\cup_{k \in \ex} S_{\mu_k(\x)}) \big) 
\]
are closed subvarieties of $Y(\B)^\circ= \MaxSpec \U(\B)$, each irreducible component of which lies in codimension $\geq 2$.
\epr

\begin{proof}
Since $\U(\B)$ is Noetherian and normal, Serre's Criterion (e.g.~\cite[Theorem 11.5.i]{Eis95}) implies that each irreducible component of the singular locus $(Y(\B)^\circ)^\sing$ has codimension $\geq2$ (this is the \emph{R1 condition}).

Because $\U(\B)$ is totally coprime, \cite[Corollary 1.7]{BFZ05} implies that the coordinate ring of the open subvariety $  S_{\x} \cup (\cup_{k \in \ex} S_{\mu_k(\x)}) \subseteq Y(\B)^\circ$ (called the \emph{upper bound cluster algebra} in \textit{loc.~cit.}) is equal to the coordinate ring of $Y(\B)^\circ$ (which is the upper cluster algebra $\U(\B)$). By the Algebraic Hartogs Lemma, the complement of $  S_{\x} \cup (\cup_{k \in \ex} S_{\mu_k(\x)})$ in $Y(\B)^\circ$ must have codimension $\geq2$; see \cite[Lemma~4.3.1]{MM15} for the general statement and \cite[Lemma~4.4.2]{MM15} for the statement in this setting.
\end{proof}
%%%%%%%%%%%%%%%%%%%%
\subsection{Proof of \thref{open-T-orb-sl}}
\label{5.5}
We prove a stronger result than that of \thref{open-T-orb-sl}, namely that $Y(\B)^\reg \backslash \VV(x)$ 
is a single torus orbit of symplectic leaves of $(Y(\B), \pi)$ for a subtorus of $T(\B)$ of rank equal to 
\[
\dim Y(\B) - 2 \rk \pi = \dim Y(\B) - \dim \Lsc,
\]
where $\Lsc$ is a symplectic leaf of $Y(\B)$ of maximal dimension.

\bth{small-T-orbit} Assume that in the setting of \thref{sections-antican}, 
$\theta = \{ i_1, \ldots, i_{N -2r} \}$ is one of the elements of $\Delta$ such that 
$\chi_\theta \neq 0$. Let $T(\B)_0 \cong (\Cset^\times)^{N-2r}$ be the subtorus 
of $T(\B) \cong (\Cset^\times)^{n(\B)}$, 
corresponding to coordinates $i_1, \ldots, i_{N-2r}$. Then the following hold:
\begin{enumerate}
\item[(i)] 
The restriction $\chi_\theta|_{(Y(\B)^\circ)^\reg}$ is a nowhere vanishing section of the 
anticanonical bundle $K^*_{(Y(\B)^\circ)^\reg}$. 
\item[(ii)] The Zariski open subset
\[
Y(\B)^\reg \backslash \VV(x) = (Y(\B)^\circ)^\reg
\]
of $Y(\B)$ is a single $T(\B)_0$-orbit of symplectic leaves of $(Y(\B), \pi)$. 
\end{enumerate}
\eth
\begin{proof} (i) The section $\chi_\theta|_{(Y(\B)^\circ)^\reg}$ does not vanish on 
\[
(Y(\B)^\circ)^\reg \backslash \big( S_{\x} \cup (\cup_{k \in \ex} S_{\mu_k(\x)}) \big)
\]
by \thref{sections-antican}. It can not vanish anywhere on $(Y(\B)^\circ)^\reg$ because otherwise its zero locus 
will be a codimension 1 subvariety of ${(Y(\B)^\circ)^\reg}$, which is disjoint from the union
$S_{\x} \cup (\cup_{k \in \ex} S_{\mu_k(\x)})$. This would contradict \prref{codim2}.

(ii) Let $\Lsc_0$ be a symplectic leaf of $(S_\x, \pi|_{S_\x})$. There exists a symplectic leaf $\Lsc$ of $(Y(\B), \pi)$ 
such that $\Lsc_0$ is a connected component of the intersection $\Lsc \cap S_\x$. The proof of \cite[Lemma 2.4]{GSV}
gives that 
\[
T(\B)_0 \cdot \Lsc_0 = S_\x.
\]
This proves the first inclusion in 
\begin{equation}
\label{2containments}
S_\x \subseteq T(\B)_0 \cdot \Lsc \subseteq (Y(\B)^\circ)^\reg.
\end{equation}
The second inclusion follows from \prref{containment}. Consider a point 
\[
y \in (Y(\B)^\circ)^\reg
\]
and let $\Lsc_y$ be a neighborhood (in the $C^\infty$ topology) of the symplectic leaf
of $(Y(\B), \pi)$ through $y$. Since
\[
w_{\nu_{i_1}, y} \wedge \ldots \wedge w_{\nu_{i_{N-2r}}, y} \wedge (\wedge^{2r} \pi_y) = \chi_{\theta, y} \neq 0, 
\]
there exists a neighborhood (in the $C^\infty$ topology) $\Osc_1$ of the identity element of the torus 
$T(\B)_0$ such that $\Osc_1 \cdot \Lsc_y$ contains a neighborhood (in the $C^\infty$ topology) of $y$ in $(Y(\B)^\circ)^\reg$. 
Since $S_\x$ is a Zariski open subset of the irreducible variety $(Y(\B)^\circ)^\reg$, 
\[
T(\B)_0 \cdot \Lsc_0  \cap  \Osc_1 \cdot \Lsc_y =
S_\x \cap \Osc_1 \cdot \Lsc_y \neq \varnothing.
\]
Hence, there exist $t \in \Osc_1$ and $t' \in T(\B)_0$ such that
\[
\psi(t) \cdot \Lsc_y \subseteq \psi(t') \cdot \Lsc. 
\]
Note that the right hand side is a symplectic leaf of $(Y(\B), \pi)$, while the left hand side is a 
$C^\infty$ open subset of a symplectic leaf of $(Y(\B), \pi)$. 
Therefore, 
\[
y \in \Lsc_y \subseteq \psi( t^{-1} t') \Lsc,
\]
which shows that the second containment in \eqref{2containments} is an equality.
\end{proof}
%%%%%%%%%%%%%%%%
\sectionnew{Poisson order structures}
\lb{qtorus}
This section contains a construction of Poisson order structures on strict root of unity upper quantum cluster algebras that will be
used to link Poisson geometry and representation theory.
%%%%%%%%%%%
\subsection{Poisson orders on mixed root of unity quantum tori}
\label{4.1}
Denote 
\[
\Zset^N_\geq := \{ (m_1, \ldots m_N) \in \Zset^N \mid m_i \geq 0, \forall i \in \ninv \}. 
\]

\bpr{Pord-tori} Assume that $\U_\ep(M_\ep, \B,\inv)$ is a strict root of unity upper quantum cluster algebra and $\ell$ is an odd positive integer that is coprime to the diagonal entries 
of the skew-symmetrizing matrix $D$ for the principal part of the exchange matrix $\B$.

Then, for every seed $(M'_\ep, \B') \sim (M_\ep, \B)$, the pair  
\[
(\TT_\ep(M_\ep)_\geq, \LL_\ep(M_\ep)_\geq)
\]
has a unique Poisson order structure $\partial' : \LL_\ep(M_\ep)_\geq \to \Der_\Cset(\TT_\ep(M_\ep)_\geq)$ such that 
\begin{equation}
\label{partial'}
\partial'_{M'_\ep(\ell f)} (M_\ep(g)) =  \frac{1}{\ell} \La_{(M'_\ep, \B')}(f,g) M_\ep(\ell f +g).
\end{equation}
The underlying Poisson structure is given by 
\begin{equation}
\label{Poisson-L2}
\{ y'_k, y'_i \}_{\partial'} = \La_{(M'_\ep, \B')}(e_k,e_i) y'_k y'_i, \quad \forall k,i \in [1,N],
\end{equation}
where $y'_k := M'_\ep(e_k)^\ell$ for $1 \leq k \leq N$. 
\epr
Recall from \prref{struct-mut} the construction of the skew-symmetric matrices $\La_{(M'_\ep, \B')} \in M_N(\Zset)$, 
identified with the corresponing bilinear forms 
\[
\La_{(M'_\ep, \B')} : \Zset^N \times \Zset^N \to \Zset.
\]
In the proposition we suppress the dependence of the $\Cset$-linear map $\partial'$ on the choice of seed $(M'_\ep, \B')$
for simplicity of the notation. Note that, in view of \eqref{Tep-geq} and \eqref{Lep-geq},
\[
\TT_\ep(M'_\ep)_\geq = \Span_\Cset \{ M'_\ep(f) \mid f \in \Zset^N_\geq\}, \quad
\LL_\ep(M'_\ep)_\geq = \Span_\Cset \{ M'_\ep(\ell f) \mid f \in \Zset^N_\geq\}. 
\]
\begin{proof} Set $\La':= \La_{(M'_\ep, \B')}$. Let $(M'_q, \B') \sim (M_q, \B)$ be the seed 
corresponding to $(M'_\ep, \B')$ under the exchange graph isomorphisms of Theorems \ref{tisom-e-graphs} and \ref{tcent-isom}(i). 
Consider the surjective $\Cset$-algebra homomorphism from \eqref{kappa} for the seed $(M'_\ep, \B')$:
\begin{equation}
\label{kappa-ep'}
\kappa'_\ep : \TT_q(M'_q)_\geq \to \TT_\ep(M'_\ep)_\geq
\end{equation}
and its $\Cset$-linear section
\begin{equation}
\label{iota}
\iota'_\ep : \TT_\ep(M'_\ep)_\geq \to \TT_q(M'_q)_\geq, \quad M'_\ep(f) \mt M'_q(f), \; \; \forall f \in \Zset^N_\geq. 
\end{equation}
By \leref{Pord-spec} and eq. \eqref{ker-kappa}, we obtain 
a Poisson order structure $\DD'$ on the pair $(\TT_\ep(M'_\ep), \ZZ(\TT_\ep(M'_\ep)))$. It is given by
\begin{align*}
\DD'_{M'_\ep(\ell f)} (M_\ep(g)) &= \kappa'_\ep \Bigg( \frac{q^{\ell \La'(f,g)/2}- q^{-\ell \La'(f,g)/2}}{q^{1/2} - \ep^{1/2}}  M'_\ep(\ell f +g) \Bigg),
\\&= c M'_\ep(\ell f +g), \quad \forall f,g \in \Zset^N_\geq,
\end{align*} 
where 
\[
c = \frac{q^{\ell \La'(f,g)/2}- q^{-\ell \La'(f,g)/2}}{q^{1/2} - \ep^{1/2}} \Bigg{|}_{q^{1/2} = \ep^{1/2}} =  2 \ell \La'(f,g) \ep^{\ell \La'(f,g)/2} \ep^{-1/2}
= 2 \ell \ep^{-1/2} \La'(f,g).
\]
The underlying Poisson structure satisfies
\[
\{ y'_k, y'_i \}_{\DD'} = 
\DD'_{M'_\ep(\ell e_k)} \big( M'_\ep(\ell e_i) \big) =
2 \ell^2 \ep^{-1/2} \La'(e_k,e_i) y'_k y'_i, \quad k,i \in [1,N].
\]
In particular, $\LL_\ep(M'_\ep)_\geq$ is a Poisson subalgebra of $\ZZ(\LL_\ep(M'_\ep)_\geq)$ with respect to the bracket 
$\{-, -\}_{\DD'}$. Therefore,
\[
\partial' :=  \frac{\ep^{1/2}}{2 \ell^2} \DD'|_{\LL_\ep(M'_\ep)_\geq}
\]
is a Poisson order structure on the pair $(\TT_\ep(M'_\ep)_\geq, \LL_\ep(M'_\ep)_\geq)$, which is given by \eqref{partial'}. 
Eq. \eqref{Poisson-L2} follows at once from \eqref{partial'}. In light of the $\Cset$-linearity of $\partial'$, there 
is a unique Poisson order structure satisfying \eqref{partial'}.
\end{proof}
%%%%%%%%%%%%
\subsection{Poisson order structures on root of unity upper quantum cluster algebras}
\label{4.2}
\bth{Uep-Pord} Assume that $\U_\ep(M_\ep, \B,\inv)$ is a strict root of unity upper quantum cluster algebra and $\ell$ is an odd positive integer that is coprime to the diagonal entries 
of the skew-symmetrizing matrix $D$ for the principal part of the exchange matrix $\B$ such that
\begin{equation}
\label{assume}
\A_\ep(M_\ep, \B,\inv) = \U_\ep(M_\ep, \B,\inv).
\end{equation}
Then the pair $(\U_\ep(M_\ep, \B,\inv), \C\U_\ep(M_\ep, \B,\inv))$ admits a Poisson order structure such that under the isomorphism 
\[
\C\U_\ep(M_\ep, \B,\inv) \cong \U(\B,\inv)
\]
from \eqref{U-CU-isom} the underlying Poisson structure on $\C\U_\ep(M_\ep, \B,\inv)$ corresponds to the GSV Poisson structure on $\U(\B,\inv)$.
\eth
\begin{proof} Recall the surjective $\Cset$-algebra homomorphism 
\[
\kappa_\ep : \TT_q(M_q)_\geq \to \TT_\ep(M_\ep)_\geq
\]
from \eqref{kappa} and consider its restriction to $\U_q(M_\ep, \B,\inv)$. For all seeds $(M'_q, \B') \sim (M_q, \B)$ and indices $k \in \ex$, we have
\begin{multline}
\label{eq-int}
(q^{1/2} - \ep^{1/2}) \TT_q(M'_q)_\geq \cap \TT_q(M''_q)_\geq  = \TT_q(M'_q)_\geq \cap (q^{1/2} - \ep^{1/2}) \TT_q(M''_q)_\geq 
\\
= (q^{1/2} - \ep^{1/2}) \big( \TT_q(M'_q)_\geq \cap \TT_q(M''_q)_\geq \big).
\end{multline}
The first equality follows from \cite[Proposition 5.10]{NTY2} and the second is an immediate consequence of the first one. 
Since $\Ker \kappa_\ep = (q^{1/2} - \ep^{1/2}) \TT_q(M_q)_\geq$ by eq. \eqref{ker-kappa},
\begin{equation}
\label{ker-kep-Upe}
\Ker  \kappa_\ep|_{\U_q(M_\ep, \B,\inv)} = (q^{1/2} - \ep^{1/2}) \U_q(M_\ep, \B,\inv). 
\end{equation}
The exchange graph isomorphism from \thref{cent-isom}(i) implies that
\begin{equation}
\label{subseteq}
\kappa_\ep \big( \U_q(M_\ep, \B,\inv) \big) \subseteq \U_\ep(M_\ep, \B,\inv).
\end{equation}
By \cite[Theorem 5.7]{NTY2}
\[
\kappa_\ep (\A_q(M_q, \B, \inv)) =\A_\ep(M_\ep, \B, \inv).
\]
Using the quantum Laurent phenomenon \eqref{qLaurent} and the assumption \eqref{assume}, we obtain
\[
\kappa_\ep (\U_q(M_q, \B, \inv)) \supseteq \kappa_\ep (\A_q(M_q, \B, \inv)) = 
\A_\ep(M_\ep, \B, \inv) = \U_\ep(M_\ep, \B, \inv).
\]
This, combined with \eqref{subseteq}, gives 
\[
\kappa_\ep (\U_q(M_q, \B, \inv)) = \U_\ep(M_\ep, \B, \inv).
\]
In view of \eqref{ker-kep-Upe} and \leref{Pord-spec}, the restriction of $\kappa_\ep$ to $\U_q(M_q, \B, \inv)$
gives rise to a Poisson order structure $\DD_\U$ on 
\[
( \U_\ep(M_q, \B, \inv), \ZZ(\U_\ep(M_q, \B, \inv)) ).
\]
Denote
\[
\partial_\U : =  \frac{\ep^{1/2}}{2 \ell^2} \DD_\U.
\]
Consider a $\Cset$-linear section 
\[
\U_\ep(M_\ep, \B, \inv) \to \U_q(M_q, \B, \inv)
\]
of the restriction of $\kappa_\ep$ to $\U_q(M_q, \B, \inv)$. For an arbitrary seed $(M'_\ep, \B') \sim (M_\ep, \B)$, 
extend it to a section of the homomorphism $\kappa'_\ep : \TT_q(M'_q)_\geq \to \TT_\ep(M'_\ep)_\geq$ from \eqref{kappa-ep'}.
\prref{Pord-tori} and the independence of the underlying Poisson structure of a Poisson order obtained by specialization (\leref{Pord-spec})
imply that 
\begin{equation}
\label{P-bra-U}
\{ y'_k, y'_i \}_{\partial_\U} =
\{ y'_k, y'_i \}_{\partial'} = \La_{(M'_\ep, \B')}(e_k,e_i) y'_k y'_i, \quad \forall k,i \in [1,N],
\end{equation}
for $y'_i := M'_\ep(e_i)^\ell$. Therefore, 
\[
\C\U_\ep(M_\ep, \B, \inv) := \bigcap_{ (M'_\ep, \B') \sim (M_\ep, \B)} \LL_\ep(M_\ep')_{\geq}
\]
is closed under the Poisson bracket $\{., .\}_{\partial_\U}$ and $\partial_U$ restricts to a Poisson order structure
\[
\partial'_\U : \C\U_\ep(M_\ep, \B, \inv) \to \Der_\Cset \big( \U_\ep(M_\ep, \B, \inv) \big).
\]
Eqs. \eqref{Poisson-L} and \eqref{P-bra-U} imply that
the underlying Poisson structure on $\C\U_\ep(M_\ep, \B,\inv)$ corresponds to the GSV Poisson structure on $\U(\B,\inv)$.
\end{proof} 
%%%%%%%%%%%%%%%%%%%%
\sectionnew{Root of unity quantum side: the fully Azumaya loci}
\lb{Azumaya}
In this section we prove our second main result describing the fully Azumaya loci of all finitely generated 
strict upper cluster algebras $\U_\ep(M_\ep, \B, \inv)$ with the property $\U_\ep = \A_\ep$. 
%%%%%%%%%%% 
\subsection{Statement of the main result}
\label{6.1}
Consider a root of unity upper quantum cluster algebra $\U_\ep(M_\ep, \B, \inv)$ and its central subalgebra
\[
\C \U_\ep(M_\ep, \B, \inv) \subset \ZZ(\U_\ep(M_\ep, \B, \inv)),
\]
given by \eqref{CU}. If, the order of $\ep$ is odd and coprime to the diagonal entries of the skew-symmetrizing matrix $D$ 
for the principal part of the exchange matrix $\B$, then by \thref{cent-isom}(iii), 
\[
\C\U_\ep(M_\ep, \B, \inv) \cong \U(\B, \inv). 
\]
If, in addition, $\U_\ep(M_\ep, \B, \inv)$  is a finitely generated $\Cset$-algebra, then by \thref{cent-isom}(iii),
it is a module finite domain over $\C\U_\ep(M_\ep, \B, \inv)$ and $\C\U_\ep(M_\ep, \B, \inv)$
is a finitely generated normal commutative $\Cset$-algebra. So, we can consider the fully Azumaya locus $\AA$
of $\U_\ep(M_\ep, \B, \inv)$ with respect to $\C\U_\ep(M_\ep, \B, \inv)$, recall \deref{fullyAzumaya}.
We will identify
\[
\MaxSpec \C\U_\ep(M_\ep, \B, \inv) \cong Y(\B) = \MaxSpec \U(\B, \inv)
\]
and think of $\AA$ as of a Zariski open subset of $Y(\B)$, recall \leref{FAzumaya-dense}. 

Recall \eqref{x-prod}. Denote the set of {\em{non-central frozen variables among the non-inverted ones}}:
\begin{equation}
\label{nc-def}
\nc := \{ i \in \ninv \mid M(e_i) \notin \ZZ(\U_\ep(M_\ep, \B, \inv)) \} 
\subset \ninv
\end{equation}
and set
\begin{equation}
\label{x-prod-nc}
x_\nc := \prod_{i \in \nc} x_i.
\end{equation}

Our main result provides an explicit characterization of the fully Azumaya locus $\AA$.
\bth{Azumaya} 
Assume that $\U_\ep(M_\ep, \B, \inv)$ is a finitely generated strict root of unity upper quantum cluster algebra such that 
the order of $\ep$ is odd and coprime to the diagonal entries of the skew-symmetrizing matrix $D$ 
for the principal part of the exchange matrix $\B$ and 
\[
\A_\ep(M_\ep, \B,\inv) = \U_\ep(M_\ep, \B,\inv).
\]
Then the fully Azumaya locus $\AA$ of 
$\U_\ep(M_\ep, \B, \inv)$ with respect to the central subalgebra
\[
\C\U_\ep(M_\ep, \B, \inv) \cong \U(\B, \inv) 
\]
satisfies
\[
(Y(\B)^\circ)^\reg \subseteq \AA \subseteq Y(\B) \backslash \VV(x_{\nc}).
\]
\eth
All important root of unity quantum cluster algebras that we are aware of do not have central frozen variables. In those situations 
\[
\nc = \ninv \quad \mbox{and} \quad
x_\nc = x,
\]
so the upper bound for the fully Azumaya locus in \thref{Azumaya} becomes
\[
Y(\B) \backslash \VV(x_{\nc}) = Y(\B) \backslash \VV(x) = Y(\B)^\circ.
\]
\bco{nc1} If $\U_\ep(M_\ep, \B, \inv)$ is a strict root of unity upper quantum cluster algebra as in 
\thref{Azumaya} such that
\[
M(e_i) \notin \ZZ(\U_\ep(M_\ep, \B, \inv)), \quad \forall i \in \ninv,
\]
then
\[
(Y(\B)^\circ)^\reg \subseteq \AA \subseteq Y(\B)^\circ.
\]
\eco
However, in general, one cannot replace the upper bound in \thref{Azumaya} 
with $Y(\B)^\circ$ as shown in \prref{nc2}.
%%%%%%%%%%% 
\subsection{Proof of the first inclusion in \thref{Azumaya}}
\label{6.2}
The first inclusion in \thref{Azumaya} follows from the following:
\bpr{Azumaya2} In the setting of \thref{Azumaya}, for all 
\[
\mm \in (Y(\B)^\circ)^\reg \subset Y(\B) \cong \MaxSpec \C\U_\ep(M_\ep, \B, \inv), 
\]
the algebra
\[
\U_\ep(M_\ep, \B, \inv) / \mm \U_\ep(M_\ep, \B, \inv)
\]
is isomorphic to 
\[
\TT_\ep(\ol{\La})/\n \TT_\ep(\ol{\La}), 
\]
where $\n$ is any maximal ideal of $\LL_\ep(\ol{\La})$ {\em{(}}recall \eqref{Tep} and \eqref{Lep}{\em{)}}
and $\La \in M_N (\Zset)$ is the skew-symmetric integer matrix compatible with $\B$ as in \deref{strict-seed}. All irreducible representations 
of the last algebra have dimension 
\begin{equation}
\label{index}
\sqrt{[\Zset^N:\Ker(\ol{\La})]},
\end{equation}
which equals the PI degree of $\U_\ep(M_\ep, \B, \inv)$.
\epr 
Here, $[G:H]$ denotes the index of a subgroups $H$ of a group $G$. 
The matrix $\ol{\La} \in M_N(\Zset/\ell)$ is identified with the corresponding bilinear 
form $\Zset^N \times \Zset^N \to \Zset/\ell$, and $\Ker (\ol{\La})$ denotes the kernel of this 
form. The index in \eqref{index} is finite because $\Ker (\ol{\La}) \supseteq (\ell \Zset)^N$. 
\begin{proof} The algebras 
\[
\U_\ep(M_\ep, \B, \inv) / \mm \U_\ep(M_\ep, \B, \inv)
\]
are isomorphic to each other for all 
\[
\mm \in (Y(\B)^\circ)^\reg.
\]
This follows by combining Theorems \ref{tUep-Pord} and \ref{topen-T-orb-sl}, \coref{P-ort-FAzumaya}, and the fact that $\varphi$ is an action of $T(\B)$ on 
$\U_\ep(M_\ep, \B, \inv)$ by algebra automorphisms that preserves the central subalgebra $\C\U_\ep(M_\ep, \B, \inv)$, see \leref{Caction}. 
So, we can restrict ourselves to the special case when 
\begin{equation}
\label{special-m}
\mm \in \MaxSpec \LL_\ep(M_\ep) \cong \MaxSpec \LL(\x) \subset (Y(\B)^\circ)^\reg.
\end{equation}
Denote by $E$ the multiplicative subset of $\C\U_\ep(M_\ep, \B, \inv)$, generated by $M_\ep(e_i)^{\ell}$ for $i \in [1,N] \backslash (\ex \sqcup \inv)$.
We have 
\begin{equation}
\label{UTT-isom0}
\U_\ep(M_\ep, \B, \inv) [E^{-1}] \cong \TT_\ep(M_\ep) \cong \TT_\ep(\ol{\La})  
\end{equation}
and, on the level of their centers,
\begin{equation}
\label{UTT-isom}
\C\U_\ep(M_\ep, \B, \inv) [E^{-1}] \cong \LL_\ep(M_\ep) \cong \LL_\ep(\ol{\La}).  
\end{equation}
For the rest of the proof assume \eqref{special-m}. This implies that $\mm$ is a maximal ideal of $\C\U_\ep(M_\ep, \B, \inv)$
that is disjoint from $E$. This property and eq. \eqref{UTT-isom0} give
\begin{align*}
&\U_\ep(M_\ep, \B, \inv) / \mm \U_\ep(M_\ep, \B, \inv) 
\\
\cong &\U_\ep(M_\ep, \B, \inv) [E^{-1}] / \mm \U_\ep(M_\ep, \B, \inv) [E^{-1}]
\cong \TT_\ep(\ol{\La})/ \n \TT_\ep(\ol{\La}). 
\end{align*}
Here $\n$ is the maximal ideal of $\LL_\ep(\ol{\La})$ corresponding to the maximal ideal 
$\mm [E^{-1}]$ of $\U_\ep(M_\ep, \B, \inv) [E^{-1}]$ under the isomorphism
\eqref{UTT-isom}.

Since $\TT_\ep(\ol{\La})$ is an Azumaya algebra of PI degree $\sqrt{[\Zset^N:\Ker(\ol{\La})]}$ by \cite[Proposition 6.1(3-4)]{HLY}, all irreducible representations 
of $\TT_\ep(\ol{\La})$ have dimension equal to the same integer. The latter equals the PI degree of $\U_\ep(M_\ep, \B, \inv)$ by \cite[Proposition 6.4]{HLY}.
\end{proof}
%%%%%%%%%%% 
\subsection{Proof of the second inclusion in \thref{Azumaya}}
\label{6.3}
Consider a strict root of unity quantum cluster algebra $\U_\ep(M_\ep, \B, \inv)$ with skew-symmetric integer matrix $\La \in M_N (\Zset)$ as in \deref{strict-seed}. 

For $j \in [1,N]\backslash (\ex \sqcup \inv)$, denote by $\La_j \in M_{N-1}(\Zset)$ the submatrix of $\La$, obtained by removing the $j$-th row and column. Set
\[
\Zset^{N-1}_j := \bigoplus_{i \neq j} \Zset e_i \subset \Zset^N.
\]
The bilinear form  $\Zset^{N-1}_j \times \Zset^{N-1}_j \to \Zset/\ell$ associated to $\ol{\La}_j$ is the restriction to $\Zset^{N-1}_j$ of the bilinear form 
$\Zset^N \times \Zset^N \to \Zset/\ell$ associated to $\ol{\La}$. 
Let $\TT_\ep(\ol{\La}_j)$ be the root of unity quantum subtorus of $\TT_\ep(\ol{\La})$ spanned by $x^f$ for $f \in \Zset^{N-1}_j$, cf. \eqref{Tep}. 
Denote
\begin{equation}
\label{J-def}
J:= \ninv \backslash \{ j \}.
\end{equation}
It is easy to see that we have a second description of $\TT_\ep(\ol{\La}_j)$:
\begin{equation}
\label{Tquot}
\big( \TT_\ep(\ol{\La})_\geq [(x^{e_i})^{-1} ; i \in J] \big) / (x^{e_j}) \cong \TT_\ep(\ol{\La}_j).
\end{equation}
Recall that an element $a$ of an algebra $R$ is called normal if $Ra = aR$. 
Here and below, for such an element $a \in R$, we denote the principal ideal 
\[
(a) := Ra = aR.
\]
\bpr{Azumaya3} Assume the setting of \thref{Azumaya}.  

Then the open subset 
\[
\big( \VV(x_j) \cap S_{\x} \big) \backslash \VV(x /x_j)
\]
of $\VV(x_j)$ is non-empty and for all 
\begin{equation}
\label{m-assume}
\mm \in \big( \VV(x_j) \cap S_{\x} \big) \backslash \VV(x /x_j) \subset Y(\B) \cong \MaxSpec \C\U_\ep(M_\ep, \B, \inv), 
\end{equation}
the algebra
\[
\U_\ep(M_\ep, \B, \inv) / \mm \U_\ep(M_\ep, \B, \inv)
\]
is isomorphic to 
\[
\TT_\ep(\ol{\La}_j)/\n \TT_\ep(\ol{\La}_j), 
\]
where $\n$ is a maximal ideal of $\LL_\ep(\ol{\La}_j)$ 
and $\La \in M_N (\Zset)$ is the skew-symmetric integer matrix as in \deref{strict-seed}. All irreducible representations 
of the last algebra have dimension 
\[
\sqrt{[\Zset^{N-1}_j:\Ker(\ol{\La}_j)]}. 
\]
\epr 
\begin{proof} Recall \eqref{J-def}.
We have 
\[
\U_\ep(M_\ep, \B, \inv) [M_\ep(e_i)^{-\ell} ; i \in J] \cong \TT_\ep(M_\ep)_{\geq} [M_\ep(e_i)^{-\ell} ; i \in J] \cong \TT_\ep(\ol{\La})_\geq [(x^{e_i})^{-\ell} ; i \in J].
\]
Combining this with \eqref{Tquot} gives
\begin{equation}
\label{UTj-isom}
\Big( \U_\ep(M_\ep, \B, \inv) [M_\ep(e_i)^{-\ell} ; i \in J] \Big)/ (M(e_j)) \cong \TT_\ep(\ol{\La}_j).
\end{equation}
Analogously, on the level of centers one shows that 
\begin{equation}
\label{isom-for-ideals}
\C\U_\ep(M_\ep, \B, \inv) [M_\ep(e_i)^{-\ell} ; i \in J] \cong \LL_\ep(\ol{\La})_\geq [(x^{e_i})^{-\ell} ; i \in J].
\end{equation}
The assumption \eqref{m-assume} implies that the maximal ideal $\mm$ of $\C\U_\ep(M_\ep, \B, \inv)$ 
is disjoint from the multiplicative set $[M_\ep(e_i)^{\ell} ; i \in J]$. 
Denote by $\n$ the maximal ideal of $\LL_\ep(\ol{\La})_\geq [x_i^{-\ell} ; i \in J]$ that corresponds to the maximal ideal
$\mm [M_\ep(e_i)^{-\ell} ; i \in J]$ of $\C\U_\ep(M_\ep, \B, \inv) [M_\ep(e_i)^{-\ell} ; i \in J]$ under the isomorphism 
\eqref{isom-for-ideals}. We have,
\begin{align*}
R:=&\big(\TT_\ep(\ol{\La})_\geq [x_i^{-\ell} ; i \in J]) / \n \big(\TT_\ep(\ol{\La})_\geq [x_i^{-\ell} ; i \in J] \big),
\\
\cong&\big( \U_\ep(M_\ep, \B, \inv) [M_\ep(e_i)^{-\ell} ; i \in J] \big) / \mm \big( \U_\ep(M_\ep, \B, \inv) [M_\ep(e_i)^{-\ell} ; i \in J] \big)
\\
\cong&\U_\ep(M_\ep, \B, \inv) / \mm \U_\ep(M_\ep, \B, \inv).
\end{align*}
It follows from \eqref{m-assume} that 
\[
(x^{e_j})^\ell \in \n, 
\]
and thus, $(x^{e_j})^\ell =0$ as an element of the algebra $R$. Therefore, all irreducible representations of $R$ are annihilated by $x^{e_j}$,
and so, those representations are in bijection with the irreducible representations of $R/(x^{e_j})$. Now we invoke \eqref{UTj-isom}, to obtain
\[
R/(x^{e_j}) \cong \big( \TT_\ep(\ol{\La})_\geq [(x^{e_i})^{-\ell} ; i \in J] \big)
/ (x^{e_j}) \n
\cong \TT_\ep(\ol{\La}_j) 
\]
for an ideal $\n'$ of $\TT_\ep(\ol{\La}_j)$. By \cite[Proposition 6.1(3-4)]{HLY}, all irreducible representations 
of $\TT_\ep(\ol{\La}_j)$ have dimension equal to 
\[
\sqrt{[\Zset^{N-1}_j:\Ker(\ol{\La}_j)]},
\]
and hence, the same holds for the algebra $\U_\ep(M_\ep, \B, \inv) / \mm \U_\ep(M_\ep, \B, \inv)$.
\end{proof}
\bpr{Azumaya4} Assume the setting of \thref{Azumaya}. Let $j \in \nc$. Let
\[
\mm \in \big( \VV(x_j) \cap S_{\x} \big) \backslash \VV(x /x_j) \subset Y(\B).
\]
The all irreducible representations of the algebra
\[
\U_\ep(M_\ep, \B, \inv) / \mm \U_\ep(M_\ep, \B, \inv)
\]
have dimension strictly less than the PI degree of $\U_\ep(M_\ep, \B, \inv)$,
\[
\sqrt{[\Zset^N:\Ker(\ol{\La})]}.
\]
\epr 
\begin{proof}
Obviously, $\Ker (\ol{\La}) \cap \Zset_j^{N-1} \subseteq \Ker ( \ol{\La}_j)$. Therefore, 
\[
[\Zset^N:\Ker(\ol{\La})] \geq [\Zset^{N-1}_j: (\Ker (\ol{\La}) \cap \Zset_j^{N-1}) ] \geq
[\Zset^{N-1}_j:\Ker(\ol{\La}_j)].
\]
The first inequality is in fact strict. If it is an equality, then $e_j \in \Ker(\ol{\La})$, which implies that 
\[
x^{e_j} \in \ZZ(\TT_\ep(\ol{\La})), 
\]
which implies that
\[
M(e_j) \in \ZZ(\U_\ep(M_\ep, \B, \inv))
\]
because $\U_\ep(M_\ep, \B, \inv)$ is a subalgebra of the skew field of fractions of $\TT_\ep(M_\ep)$. 
This contradicts the definition \eqref{nc-def} of the set $\nc$. Hence, 
\[
\sqrt{[\Zset^N:\Ker(\ol{\La}_j)]} <
\sqrt{[\Zset^N:\Ker(\ol{\La})]}.
\]
and the proposition now follows from \prref{Azumaya3}.
\end{proof}
With this proposition we complete the proof of \thref{Azumaya}:
\medskip

\noindent
{\em{Proof of the second inclusion in \thref{Azumaya}.}} Let $j \in \nc$. Propositions \ref{pAzumaya3} and 
\ref{pAzumaya4} imply that
\[
\big( \VV(x_j) \cap S_{\x} \big) \backslash \VV(x /x_j) \subseteq Y(\B) \backslash \AA
\]
and that the first set is a non-empty Zariski open subset of $\VV(\x_j)$. 
Since $\VV(x_j)$ is irreducible by \prref{p-primes}, 
\[
\ol{\big( \VV(x_j) \cap S_{\x} \big) \backslash \VV(x /x_j)} = \VV(x_j).
\]
Furthermore, \leref{FAzumaya-dense} implies that $Y(\B) \backslash \AA$ 
is a Zariski closed subset of $Y(\B)$. Therefore,
\[
\VV(x_j) \subseteq Y(\B) \backslash \AA, \quad \forall j \in \nc.
\]
hence, $\VV(x_\nc) \subseteq Y(\B) \backslash \AA$, and thus 
$\AA \subseteq Y(\B) \backslash \VV(x_{\nc})$.
\qed
\medskip

Last we prove that one cannot replace the upper bound in \thref{Azumaya} 
with $Y(\B)^\circ$. 
\bpr{nc2} Assume the setting of \thref{Azumaya}. Let
\[
j \in \ninv \backslash \nc.
\]
All points 
\[
\mm \in \Big( \VV(x_j) \cap \Big( \bigcup_{(\x', \B') \sim (\x,\B)} S_{\x'} \Big) \Big) \backslash \VV(x/x_j)
\]
belong to the fully Azumaya locus of $\U_\ep(M_\ep, \B, \inv)$ with respect to 
$\C\U_\ep(M_\ep, \B, \inv)$.
\epr
\begin{proof}
It is sufficient to prove the statement for 
\[
\mm \in (\VV(x_j) \cap S_{\x} ) \backslash \VV(x/x_j).
\]
Applying \prref{Azumaya3}, we obtain that all irreducible representations 
of the algebra $\U_\ep(M_\ep, \B, \inv) /\mm \U_\ep(M_\ep, \B, \inv)$
have dimension 
\[
\sqrt{[\Zset^{N-1}_j:\Ker(\ol{\La}_j)]},
\]
Since $j \in \ninv \backslash \nc$, we have $M(e_j) \in \ZZ(U_\ep(M_\ep, \B, \inv))$, and thus, 
\[
e_j \in \Ker (\ol{\La}). 
\]
Therefore, 
\[
[\Zset^N:\Ker(\ol{\La})] = [\Zset^{N-1}_j: (\Ker (\ol{\La}) \cap \Zset_j^{N-1}) ] =
[\Zset^{N-1}_j:\Ker(\ol{\La}_j)].
\]
Hence, 
\[
\sqrt{[\Zset^N:\Ker(\ol{\La}_j)]} <
\sqrt{[\Zset^N:\Ker(\ol{\La})]}.
\]
and the proposition now follows from \prref{Azumaya3}.
\end{proof}
%%%%%%%%%%%%%%%%%%%%%%%%%%%%%%%%
\sectionnew{Special cases and examples}
\label{spec-cases}
%%%%%%%%%%
\subsection{The Richardson divisor of a Schubert cell}
\label{7.1}
The explicit Zariski open torus orbit of symplectic leaves from \thref{open-T-orb-sl} is a far reaching generalization of the 
complement of the Richardson divisor of a Schubert cell for a complex simple Lie group $G$. 

Let $G$ be a complex simple Lie group with a pair of opposite Borel subgroups $B_\pm$. Denote by 
$T := B_+ \cap B_-$ the corresponding maximal torus of $G$ and by $U_\pm$ the unipotent radicals of $B_\pm$. 
Let $W$ be the Weyl group of $G$, identified with $N(T)/T$, where $N(T)$ is the normalizer of $T$ in $G$. 
Denote by $s_1, \ldots, s_r$ the set of simple reflections of $W$ and by $l : W \to \Zset_{\geq 0}$ 
the length function on $W$. 

The full flag variety has a canonical Poisson structure $\pi$ which is the descent of the standard Poison--Lie 
structure on $G$, \cite[Sect. 1.2-1.3]{CP}. This Poisson structure is invariant under $T$.

The {\em{Schubert cell}} of the full flag variety $G/B_+$ corresponding to $w \in W$ is 
\[
X_w^\circ := B_+ w B_+/B_+ \subset G/B_+.
\]
The {\em{open Richardson variety}} \cite{KLS} in $G/B_+$ corresponding to the pair $(w,u) \in W \times W$ is
\[
R_{w,u} := \big( B_+ w B_+ \cap B_- u B_+ \big) /B_+ \subset G/B_+.
\]
It is non-empty precisely when $u \leq w$ in the Bruhat order. We have the partition:
\[
X_w^\circ = \bigsqcup_{u \in W, u \leq W} R_{w,u}. 
\]
Denote the support of $w \in W$,
\[
\Supp w := \{ i \in [1,r] \mid s_i \leq w\}. 
\]
This is precisely the set of indices $i \in [1,r]$ such that $s_i$ appears in one, and thus, in any reduced decomposition of $w$.
The {\em{Richardson divisor}} of $X_w^\circ$ is 
\[
RD_w : = \bigcup_{i \in \Supp(w)} \Cl_{X_w^\circ}(R_{w,s_i}) =  \bigsqcup_{u \in W, u \leq w, u \neq 1} R_{w,u},
\]
where $\Cl_Y(Z)$ stands for the Zariski closure of $Z$ in $Y$. 

\bth{T-orbits} (Goodearl--Yakimov) \cite[Theorem 0.4]{GY0} The $T$-orbits of symplectic leaves of the Schubert cell $(X_w^\circ, \pi_{X_w^\circ})$ are 
the open Richardson varieties $R_{w,u}$ for $u \in W$, $u \leq w$. There is a Zariski open $T$-orbit of leaves, which is the complement 
of the Richardson divisor of the Schubert cell $X_w^\circ$:
\[
R_{w,1} = X_w^\circ \backslash RD_w.
\]
\eth
\bex{Richardson} The complement  of the Richardson divisor $RD_w$ of the Schubert cell $X_w^\circ$ is a very special case of the Zariski open torus orbit 
of symplectic leaves of $Y(\B) = \MaxSpec U(\B, \inv)$ from \thref{open-T-orb-sl}, as we show next. 
(This can be shown in the more general case of symmetrizabe Kac--Moody groups, but the setting requires more technical details.)

To each $w \in W$, one associates an exchange matrix $\B_w$ of size 
\[
l(w) \times (l(w) - |\Supp(w)|)
\]
by the first display before Theorem 10.1 of \cite{GY1}. It is known \cite{GLS1,GY3,SW} that for the corresponding cluster algebras 
without inverted frozen variables (i.e., $\inv = \varnothing$):
\[
\Cset[X_w^\circ] = \A(\B_w, \varnothing) = \U(\B_w, \varnothing).
\]
Thus,
\begin{equation}
\label{YX}
Y(\B_w) \cong X_w^\circ, 
\end{equation}
and in particular,
\[
Y(\B_w)^\sing = \varnothing, \quad Y(\B_w)^\reg = X_w^\sing.
\]

The frozen variables of the cluster algebra $\A(\B_w, \varnothing)$ are the generalized minors (cf. \cite[Eq. (2.5)]{BFZ05}),
\[
\Delta_{\vpi, w \vpi}, \quad i \in \Supp(w),
\]
considered as functions on $X_w^\circ$ via the identification $U_+ \cap w U_- w^{-1} \cong B_+ w B_+/B_+$, 
given by $g \mt g. wB_+/B_+$. There is a canonical isomorphism $T(\B_w) \cong T$ and the corresponding 
actions on \eqref{YX} coincide.

The exchange matrix $\B_w$ admits a compatible 
skew-symmetric matrix $\La_w \in M_{l(w)}(\Zset)$ and the cluster algebra $\A(\B_w, \varnothing)$,
admits a quantization, isomorphic to the corresponding integral quantum 
unipotent cell \cite{GY2}. By specialization, the GSV Poisson structure 
on $Y(\B_w) = \MaxSpec \U(\B_w, \varnothing) \cong X_w^\circ$ associated to $\La_w$ 
concides with $\pi$. We have
\[
\VV(\Delta_{\vpi, w \vpi}) = R_{w,s_i}, \quad \forall i \in \Supp(w).
\]
Thus,
\[
X_w^\circ \backslash RD_w = Y(\B)^\circ = (Y(\B)^\circ)^\reg,
\]
which shows how the Zariski open $T$-orbit of symplectic leaves in \thref{T-orbits} is a special case of that in 
\thref{open-T-orb-sl} for a small class of cluster algebras.
\eex
%%%%%%%%%%
\subsection{All frozen variables inverted}
\label{7.2}
There is an important special case of Theorems \ref{topen-T-orb-sl} and \ref{tAzumaya} when all frozen variables are inverted. 
In those situations our results take on a particularly strong form. On the Poisson side we have:
\bco{open-T-orb-sl} If $\U(\B)$ is a finitely generated upper cluster algebra with all frozen variables inverted for which there exists a skew-symmetric matrix $\La \in M_N(\Zset)$ 
such that $(\La, \B)$ is a compatible pair, then the non-singular part of $\MaxSpec U(\B)$ is a single $T(\B)$-orbit of symplectic leaves of 
the affine Poisson variety $(\MaxSpec U(\B), \pi)$.
\eco

On the Azumaya loci side we have:
\bco{Azumaya-inv} 
If $\U_\ep(M_\ep, \B)$ is a finitely generated strict root of unity upper quantum cluster algebra 
with all frozen variables inverted such that 
the order of $\ep$ is odd and coprime to the diagonal entries of the skew-symmetrizing matrix $D$ 
for the principal part of the exchange matrix $\B$ and $\U_\ep(M_\ep, \B) = \A_\ep(M_\ep, \B)$,
then fully Azumaya locus of $\U_\ep(M_\ep, \B)$
over $\C\U_\ep(M_\ep,\B)$ contains the nonsingular part of $\MaxSpec \U(\B)$.
\eco
%%%%%%%%%%
\subsection{Acyclic cluster algebras}
\label{7.3}
Another important special case of Theorems \ref{topen-T-orb-sl} and \ref{tAzumaya} is the case of acyclic cluster algebras 
when  $\U_\ep(M_\ep,\B)=\A_\ep(M_\ep,\B)$ and
an explicit presentation of $\U(\B)$ and $\U_\ep(M_\ep, \B)$ can be given, which in particular implies that those algebras
are finitely generated.

Recall that the sign pattern of an exchange matrix $\B$ is encoded in the graph $\Ga(\B)$ with vertex set $\ex$ and directed edges $(i,j)$ 
for the vertices $i,j$ with $b_{ij} >0$.  We say that $\B$ (and the corresponding cluster algebras of various kinds) are acyclic 
if $\Ga(\B)$ has no oriented cycles, cf. \cite[Definition 1.14]{BFZ05}.

\bth{acycl1} (Berenstein-Fomin-Zelevinsky) \cite[Theorems 1.18 and 1.20]{BFZ05} If $\B$ is an acyclic exchange matrix
and all variables are exchangeable {\em{(}}$\ex = [1,N]${\em{)}}, then $\A(\B)=\U(\B)$ and this algebra is isomorphic to the $\Cset$-algebra 
with generators $x_1, x'_1, \ldots, x_n, x'_n$ and relations
\begin{equation}
\label{x'x-rel}
x'_k x_k = \prod_{i, b_{ik}>0} x_i^{b_{ik}} + \prod_{i, b_{ik}<0} x_i^{-b_{ik}}, \quad \forall k \in [1,N].
\end{equation}
\eth

Analogously to this result and \cite[Theorem 7.5]{BerZe}, one proves the following:

\bpr{acycl2} If $\A_\ep(M_\ep, \B)$ is an acyclic root of unity quantum cluster algebra for which
all variables are exchangeable {\em{(}}$\ex = [1,N]${\em{)}}, then $\A_\ep(M_\ep,\B)=\U_\ep(M_\ep,\B)$ and this algebra 
is isomorphic to the $\Cset$-algebra with generators $y_1, y'_1, \ldots, y_n, y'_n$ and relations
\begin{align}
y_j y_k &= \ep^{\la_{jk}} y_k y_j, \quad &&\forall 1 \leq j < k \leq N,
\label{yy-rel} 
\\
y'_k y_k &= \ep^{\mu_{ki}/2} \prod_{i, b_{ik}>0} y_i^{b_{ik}} + \ep^{\nu_{ki}/2} \prod_{i, b_{ik}<0} y_i^{-b_{ik}}, \quad &&\forall k \in [1,N],
\label{y'y-rel}
\end{align}
where
\begin{align*}
\mu_{ki} &= \sum_{i<j, \, b_{ik}>0, b_{jk} >0} b_{ik} b_{jk} \phi_{ij} - \sum_{i, b_{ik}>0 } b_{ik} \phi_{ik},
\\
\nu_{ki} &= \sum_{i<j, \, b_{ik}<0, b_{jk} <0} b_{ik} b_{jk} \phi_{ij} + \sum_{i, b_{ik}<0 } b_{ik} \phi_{ik}
\end{align*}
and $(\phi_{ij}) \in M_N(\Zset/\ell)$ is the matrix of the root of unity toric frame $M_\ep$. 
\epr
An immediate consequence of \thref{open-T-orb-sl} is the following:

\bco{Poisson-acycl} If an acyclic exchange matrix $\B$ has a compatible skew-symmetric matrix $\La \in M_N(\Zset)$  
and all variables are exchangeable {\em{(}}$\ex = [1,N]${\em{)}}, then the maximum spectrum $Y(\B)$ 
of the algebra with generators $x_1, x'_1, \ldots, x_n, x'_n$ and relations \eqref{x'x-rel} is an affine 
Poisson manifold whose non-singular part $Y(\B)^\reg$ is a single $T(\B)$-orbit of symplectic leaves. 
\eco
An immediate consequence of \thref{Azumaya} is the following:

\bco{Azumaya-inv2} 
Assume that $\B$ is an acyclic exchange matrix, $(\la_{ij}) \in M_N(\Zset)$ is a compatible 
skew-symmetric matrix and $\ep$ is an odd root of unity whose order is coprime to the diagonal entries 
of the skew-symmetrizing matrix $D$ for the principal part of $\B$. Consider the $\Cset$-algebra 
$\U_\ep$ with generators $y_1, y'_1, \ldots, y_n, y'_n$ and relations \eqref{yy-rel}--\eqref{y'y-rel} with $\phi_{ij}:= \ol{\la}_{ij}$. 
Then the following hold:
\begin{itemize} 
\item[(i)] The unital subalgebra $\C\U_\ep$ of $\U_\ep$ generated by $y_1^\ell, (y'_1)^\ell, \ldots, y_n^\ell, (y'_n)^\ell$ is central and isomorphic 
to the algebra with generators $x_1, x'_1, \ldots, x_n, x'_n$ and relations \eqref{x'x-rel} 
via $y_k^\ell \mt x_k$, $(y'_k)^\ell \mt x_k$. $\U_\ep$ is a finitely generated module over $\C\U_\ep$. 
\item[(ii)] The fully Azumaya locus of $\U_\ep$ with respect to $\C\U_\ep$ contains 
the non-singular part of $\MaxSpec(\C\U_\ep)$.
\end{itemize}
\eco
%%%%%%%%%%
\subsection{An acyclic example}
As an example of an acyclic cluster algebra, consider the following compatible pair.
\[ \B = B=
\begin{bmatrix}
0 & -2 \\
2 & 0
\end{bmatrix}
\hspace{1cm}
\Lambda =
\begin{bmatrix}
0 & 1 \\
-1 & 0
\end{bmatrix}
\]
As in the previous subsection, every cluster variable is mutable, and so, $\ex = [1,N] = [1,2]$ and $\inv=\varnothing$. 
Because $\Ker(\B^\top)=0$, the torus $T(\B)$ is trivial, 
and so \thref{open-T-orb-sl} predicts that $Y(\B)^{\reg}$ is a single symplectic leaf.

The quiver of $B$ is the Kronecker quiver, which is acyclic. By \thref{acycl1},
\[ 
\A(\B) = \U(\B)
\cong \Cset[x_1,x_2,x_1',x_2']/( x_1x_1'-x_2^2-1,x_2x_2' - x_1^2-1).
\]
There is an even nicer presentation in terms of the following element
\[ 
z := x_1'x_2' - x_1x_2 \in \U(\B). 
\]
This element satisfies several notable identities.
\[
x_1x_2 z = x_1^2+x_2^2+1,
\hspace{1cm}
x_1 z = x_2+x_2',
\hspace{1cm}
x_2 z = x_1+x_1'.
\]
The latter two equations imply that $x_1,x_2$, and $z$ generate $\U(\B)$, and the first equation implies that the relations 
among these generators are generated by a single element; therefore,
\begin{equation}
\label{AUrel}
\A(\B) = \U(\B) \cong \Cset[x_1,x_2,z]/( x_1x_2 z - (x_1^2 + x_2^2 + 1) ).
\end{equation}
Since  
\[
f(x_1,x_2,z) := x_1x_2z - (x_1^2 + x_2^2 + 1)
\]
is an irreducible polynomial, the maximum spectrum of $Y(\widetilde{B})$ is isomorphic to the irreducible algebraic variety 
$\VV(x_1x_2 z-(x_1^2+x_2^2+1)) \subset \Cset^3$ defined by $f(x_1,x_2,z)=0$.

The GSV Poisson structure on $Y(\widetilde{B})$ extends to the Poisson structure on $\mathbb{C}^3$ with potential $f(x_1, x_2, z)$. 
The latter is defined on the coordinates by
\[ 
\{x_1,x_2\} := f_z= x_1x_2
,\;\;\;
\{x_1,z\} := - f_{x_2} = 2x_2-x_1 z
,\;\;\;
\{x_2, z\} := f_{x_1} = x_2 z - 2x_1.
\]
At each point in $Y(\B)$, the latter two functions cannot simultaneously be zero, since
\[
%\{x_1^2-x_2^2,z\} =
x_1\{x_1,z\} -x_2\{x_2,z\} =
x_1(x_2 z-2x_1) + x_2(x_1 z-2x_1)
=2x_1x_2 z -2x_1^2 -2x_2^2
\]
equals $2$ everywhere on $Y(\widetilde{B})$. This has two important consequences.
\begin{itemize}
    \item The differential of the polynomial $f(x_1,x_2, z)$ is
    \[ df = \{x_2, z\} dx_1 +\{z,x_1\}dx_2 + \{x_1,x_2\} d z. \]
    Since $df$ does not vanish on $\VV(f)$, it is a smooth variety, and so $Y(\B)^{\reg}=Y(\B)$.
    \item The Hamiltonian vector field of $z$ is
    \[
    % H_{x_1}
    % = \langle 0, x_1x_2, 2x_2-x_1 z\rangle
    % ,\;\;\;
    % H_{x_2}
    % = \langle -x_1x_2, 0, x_2 z-2x_1 \rangle
    % ,\;\;\;
    H_{z}
    = \langle \{z,x_1\}, \{z,x_2\},0\rangle.
%    = \langle x_1 z - 2x_2 , 2x_1  -x_2 z, 0 \rangle
    \]
   Since $H_z$ does not vanish on $Y(\B)$, the rank of the Poisson bracket cannot be 0 anywhere on $Y(\B)$. Since  $Y(\B)$ is smooth and 2-dimensional, the rank must be $2$ everywhere.
\end{itemize}
Since $Y(\B)$ is connected, it is a single symplectic leaf as predicted.

An interesting feature of this example is that the cluster tori do not cover $Y(\B)=Y(\B)^{\reg}$.
%To characterize the deep points of $Y(\widetilde{B})$, we first describe the clusters. 
By \cite{SZ04}, the cluster variables may be indexed by the integers so that the clusters are pairs of adjacent variables 
$(x_i,x_{i+1})$. These cluster variables may be defined recursively by the mutation identity
\[
x_{n-1}x_{n+1} = x_n^2+1.
% \text{ and }
% x_i z = x_{i-1} + x_{i+1}
\]
% This identity implies that at most one of $x_n$ and $x_{n+1}$ can vanish at the same point in $Y(\widetilde{B})$. However, at any deep point in $Y(\widetilde{B})$, at least one of $x_i$ and $x_{i+1}$ must vanish by definition. Therefore, at any deep point, either $x_{2n}$ vanishes for all $n\in \Zset$ or $x_{2n+1}$ vanishes for all $n\in \Zset$.
%
If $x_{n}=0$ at a point, then the mutation identities force $x_{n-1}x_{n+1}=1$ and
$x_{n-1}^2+1=0$, and so $x_{n-1}=\pm i$ and $x_{n+1}=\mp i$. It follows that there are four points in 
$Y(\B)$ which are not in any cluster torus, on which the sequence of cluster variables takes a periodic sequence of values
\[...0,i,0,-i,0,-i,0,i,0,...\]
Under the embedding $(x_1,x_2,z):Y(\B) \hra \Cset^3$, these points are sent to
\[ (0,\pm i,0), (\pm i,0,0) \in \Cset^3. \]

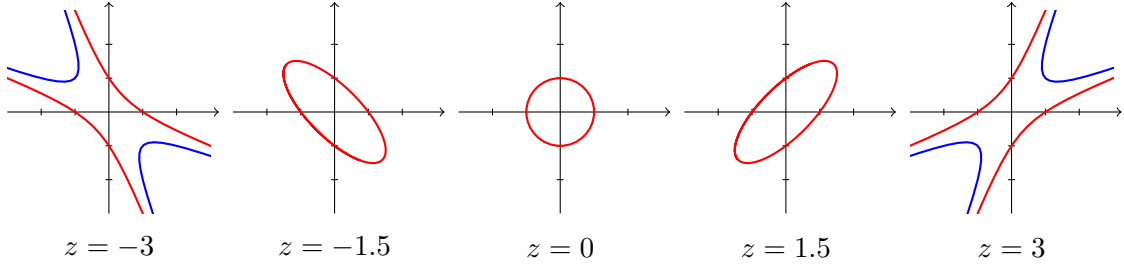
\begin{figure}[h!t]
\begin{centering}
\begin{tikzpicture}
    \begin{scope}[xshift=-6cm,scale=.45]
        \draw[->] (-3,0) to (3.25,0);
        \draw[->] (0,-3) to (0,3.25);
        \draw (-2,-.1) to (-2,.1);
        \draw (-1,-.1) to (-1,.1);
        \draw (1,-.1) to (1,.1);
        \draw (2,-.1) to (2,.1);
        \draw (-.1,-2) to (.1,-2);
        \draw (-.1,-1) to (.1,-1);
        \draw (-.1,1) to (.1,1);
        \draw (-.1,2) to (.1,2);
        \node at (0,-4) {$z=-3$};
        \clip (-3,-3) rectangle (3,3);
        \draw[thick,blue,variable=\t,domain=-3:2] plot ({2*cosh(\t)/sqrt(5)},{(-3*cosh(\t)+sqrt(5)*sinh(\t))/sqrt(5)});
        \draw[thick,blue,variable=\t,domain=-3:2] plot ({-2*cosh(\t)/sqrt(5)},{(3*cosh(\t)+sqrt(5)*sinh(\t))/sqrt(5)});
        \draw[thick,red,variable=\t,domain=-3:2] plot ({-2*sinh(\t)/sqrt(5)},{(3*sinh(\t)+sqrt(5)*cosh(\t))/sqrt(5)});
        \draw[thick,red,variable=\t,domain=-3:2] plot ({2*sinh(\t)/sqrt(5)},{-(3*sinh(\t)+sqrt(5)*cosh(\t))/sqrt(5)});
    \end{scope}
    \begin{scope}[xshift=-3cm,scale=.45]
        \draw[->] (-3,0) to (3.25,0);
        \draw[->] (0,-3) to (0,3.25);
        \draw (-2,-.1) to (-2,.1);
        \draw (-1,-.1) to (-1,.1);
        \draw (1,-.1) to (1,.1);
        \draw (2,-.1) to (2,.1);
        \draw (-.1,-2) to (.1,-2);
        \draw (-.1,-1) to (.1,-1);
        \draw (-.1,1) to (.1,1);
        \draw (-.1,2) to (.1,2);
        \node at (0,-4) {$z=-1.5$};
        \path[use as bounding box] (-3,-3) rectangle (3,3);
        \draw[thick,red,variable=\t,domain=-300:200,samples=100] plot ({2*cos(\t)/sqrt(4-1.5^2)},{-(1.5*cos(\t)+sqrt(4-1.5^2)*sin(\t))/(sqrt(4-1.5^2))});
    \end{scope}
    \begin{scope}[xshift=0cm,scale=.45]
        \draw[->] (-3,0) to (3.25,0);
        \draw[->] (0,-3) to (0,3.25);
        \draw (-2,-.1) to (-2,.1);
        \draw (-1,-.1) to (-1,.1);
        \draw (1,-.1) to (1,.1);
        \draw (2,-.1) to (2,.1);
        \draw (-.1,-2) to (.1,-2);
        \draw (-.1,-1) to (.1,-1);
        \draw (-.1,1) to (.1,1);
        \draw (-.1,2) to (.1,2);
        \node at (0,-4) {$z=0$};
        \clip (-3,-3) rectangle (3,3);
        \draw[thick,red] (0,0) circle (1);
    \end{scope}
    \begin{scope}[xshift=3cm,scale=.45]
        \draw[->] (-3,0) to (3.25,0);
        \draw[->] (0,-3) to (0,3.25);
        \draw (-2,-.1) to (-2,.1);
        \draw (-1,-.1) to (-1,.1);
        \draw (1,-.1) to (1,.1);
        \draw (2,-.1) to (2,.1);
        \draw (-.1,-2) to (.1,-2);
        \draw (-.1,-1) to (.1,-1);
        \draw (-.1,1) to (.1,1);
        \draw (-.1,2) to (.1,2);
        \node at (0,-4) {$z=1.5$};
        \path[use as bounding box] (-3,-3) rectangle (3,3);
        \draw[thick,red,variable=\t,domain=-300:200,samples=100] plot ({2*cos(\t)/sqrt(4-1.5^2)},{(1.5*cos(\t)+sqrt(4-1.5^2)*sin(\t))/(sqrt(4-1.5^2))});
    \end{scope}
    \begin{scope}[xshift=6cm,scale=.45]
        \draw[->] (-3,0) to (3.25,0);
        \draw[->] (0,-3) to (0,3.25);
        \draw (-2,-.1) to (-2,.1);
        \draw (-1,-.1) to (-1,.1);
        \draw (1,-.1) to (1,.1);
        \draw (2,-.1) to (2,.1);
        \draw (-.1,-2) to (.1,-2);
        \draw (-.1,-1) to (.1,-1);
        \draw (-.1,1) to (.1,1);
        \draw (-.1,2) to (.1,2);
        \node at (0,-4) {$z=3$};
        \clip (-3,-3) rectangle (3,3);
        \draw[thick,blue,variable=\t,domain=-3:2] plot ({2*cosh(\t)/sqrt(5)},{(3*cosh(\t)+sqrt(5)*sinh(\t))/sqrt(5)});
        \draw[thick,blue,variable=\t,domain=-3:2] plot ({-2*cosh(\t)/sqrt(5)},{(-3*cosh(\t)+sqrt(5)*sinh(\t))/sqrt(5)});
        \draw[thick,red,variable=\t,domain=-3:2] plot ({2*sinh(\t)/sqrt(5)},{(3*sinh(\t)+sqrt(5)*cosh(\t))/sqrt(5)});
        \draw[thick,red,variable=\t,domain=-3:2] plot ({-2*sinh(\t)/sqrt(5)},{-(3*sinh(\t)+sqrt(5)*cosh(\t))/sqrt(5)});
        % \draw[thick,blue,variable=\t,domain=-3:2] plot ({(e^(\t)+e^(-\t))/sqrt(5)},{((3+sqrt(5))*e^(\t)+(3-sqrt(5))*e^(-\t))/(2*sqrt(5))});
        % \draw[thick,blue,variable=\t,domain=-3:2] plot ({-(e^(\t)+e^(-\t))/sqrt(5)},{-((3+sqrt(5))*e^(\t)+(3-sqrt(5))*e^(-\t))/(2*sqrt(5))});
    \end{scope}
\end{tikzpicture}
\end{centering}
\caption{Five curves in the family $\VV(x_1^2- z x_1x_2 +x_2^2+1)$ (real part in blue and imaginary part in red)}
\label{fig: pencilofconics}
\end{figure}

The embedding $Y(\B)\hra\Cset^3$ and these four points can be visualized as follows. 
Fixing a value of $z$ is equivalent to intersecting the image of $Y(\B)$ with a plane, and the result is the curve
\[\VV(x_1^2- z x_1x_2 +x_2^2+1)\subset \Cset^2.\]%, in which the fiber over a fixed value of $z$ is $\mathcal{V}(x_1^2- z x_1x_2 +x_2^2+1)$.
For every value of $z$, this curve is a conic that passes through the four points $(0,\pm i)$ and $(\pm i,0)$; Figure \ref{fig: pencilofconics} depicts these curves for five values of $z$. One may show that every conic through the points $(0,\pm i)$ and $(\pm i,0)$ appears as $z$ varies except one: the singular conic $\VV(x_1x_2)$.\footnote{There are two other singular conics through the four points, which correspond to $z=2$ and $z=-2$.}

For historical reasons, this family of curves is called the \textbf{pencil of conics} through the four \textbf{base points}  $(0,\pm i)$ and $(\pm i,0)$.
The cluster $Y(\B)$ is then identified with the total space of the pencil of conics through $(0,\pm i)$ and $(\pm i,0)$ minus one of the three singular fibers. 
Under this identification, the four points not in any cluster torus correspond to the base points in the $z=0$ fiber, which is the complex circle of radius -1.

Let $\ep$ be a primitive root of unity of odd order. \prref{acycl2} implies that 
for the root of unity quantum cluster algebra and root of unity upper quantum cluster algebra 
associated to $\B$ we have 
\begin{equation}
\A_\ep(M_\ep,\B) = \U_\ep(M_\ep,\B) \cong
\frac{\Cset \langle x_1,x_2,x_1',x_2' \rangle}{(x_1 x_2 - \ep x_2 x_1, x_1x_1'- \ep^{-1} x_2^2-1, x_2x_2' - \ep x_1^2-1)} \cdot
\label{AUrel-ep}
\end{equation}
By \coref{Azumaya-inv2} and the first part of the example, the subalgebra of this algebra generated by $y_1^\ell, (y'_1)^\ell, y_2^\ell, (y'_2)^\ell$ 
is central and is isomorphic to the algebra \eqref{AUrel} via $y_k^\ell \mt x_k$, $(y'_k)^\ell \mt x'_k$. Furthermore, 
all irreducible representations of the algebra \eqref{AUrel-ep} have dimension $\ell$, and thus, this algebra is Azumaya.
%%%%%%%%%%%%%%%%%%%%%%%%%%%%%%%%

\end{document}